\documentclass{amsart}

\usepackage{times,epsf}

\begin{document}

\newtheorem{thm}{Theorem}[section]
\newtheorem{lem}[thm]{Lemma}
\newtheorem{cor}[thm]{Corollary}

\theoremstyle{definition}
\newtheorem{defn}{Definition}[section]

\theoremstyle{remark}
\newtheorem{rmk}{Remark}[section]

\def\square{\hfill${\vcenter{\vbox{\hrule height.4pt \hbox{\vrule
width.4pt height7pt \kern7pt \vrule width.4pt} \hrule height.4pt}}}$}
\def\T{\mathcal T}

\newenvironment{pf}{{\it Proof:}\quad}{\square \vskip 12pt}

\title{Uniform 1-cochains and Genuine Laminations}
\author{Baris Coskunuzer}
\address{Department of Mathematics \\ Princeton University \\ Princeton, NJ 08544}
\email{baris@math.princeton.edu}

\maketitle

%% User definitions:

\newcommand{\Angle}{\operatorname{angle}}
\newcommand{\Radius}{\operatorname{radius}}
\newcommand{\Sup}{\operatorname{Sup}}
\newcommand{\Int}{\operatorname{int}}
\newcommand{\Length}{\operatorname{length}}
\newcommand{\Area}{\operatorname{Area}}
\newcommand{\Lim}{\operatorname{Lim}}
\newcommand{\Out}{\operatorname{Out}}
\newcommand{\Isom}{\operatorname{Isom}}
\newcommand{\Homeo}{\operatorname{Homeo}}
\newcommand{\ID}{\operatorname{id}}
\newcommand{\cirB}{\overset{\circ}{B}}
\newcommand{\cirC}{\overset{\circ}{C}}
\newcommand{\cirD}{\overset{\circ}{D}}
\newcommand{\cirE}{\overset{\circ}{E}}
\newcommand{\cirI}{\overset{\circ}{I}}
\newcommand{\cirN}{\overset{\circ}{N}}
\newcommand{\cirT}{\overset{\circ}{T}}
\newcommand{\cirV}{\overset{\circ}{V}}
\newcommand{\cirW}{\overset{\circ}{W}}
\newcommand{\cirY}{\overset{\circ}{Y}}
\newcommand{\cirtau}{\overset{\circ}{\tau} }
\newcommand{\Si}{S^2_{\infty }}
\newcommand{\kappaj}{\kappa _{j+1}}
%%  argument form of \kappaj in amstex source: +1
\newcommand{\lambdaij}{\lambda _{ij}}
\newcommand{\sigmaij}{\sigma _{ij}}
\newcommand{\Dij}{D_{ij}}
\newcommand{\Eij}{E_{ij}}
\newcommand{\Hij}{H_{ij}}
\newcommand{\Yj}{Y_{j+1}}
%%  argument form of \Yj in amstex source: +1
\newcommand{\solidtorus}{D^2\times S^1}
\newcommand{\finv}{f^{-1}}
\newcommand{\kinv}{k^{-1}}
\newcommand{\hinv}{h^{-1}}
\newcommand{\ginv}{g^{-1}}
\newcommand{\gprimeinv}{g^{\prime -1}}
\newcommand{\pinv}{p^{-1}}
\newcommand{\qinv}{q^{-1}}
\renewcommand{\Bbb}{\mathbf}
\newcommand{\hyp}{\Bbb H}
\newcommand{\BR}{\Bbb R}
\newcommand{\BB}{\Bbb B}
\newcommand{\BN}{\Bbb N}
\newcommand{\BZ}{\Bbb Z}
\newcommand{\osolidtorus}{\cirD {}^2\times S^1}
\newcommand{\Vprime}{V^\prime }
\newcommand{\Viprime}{V_i^\prime }
\newcommand{\Vidprime}{V_i^{\prime \prime }}
\newcommand{\Vdprime}{V^{\prime \prime }}
\newcommand{\sigmaijprime}{\sigma _{ij}^\prime }
\newcommand{\Hijprime}{H_{ij}^\prime }
\newcommand{\Hijdprime}{H_{ij}^{\prime \prime }}
\newcommand{\Hiprime}{H_i^\prime }
\newcommand{\Hidprime}{H_i^{\prime \prime }}

\begin{abstract}

We construct a pair of transverse genuine laminations on an atoroidal 3-manifold
admitting transversely orientable uniform 1-cochain. The laminations are induced by the uniform 1-cochain and
they are indeed the "straightening" of the coarse laminations defined in [Ca], by using 
minimal surface techniques.
Moreover, when you collapse these laminations, you can get a topological pseudo-Anosov 
flow, as defined by Mosher, [Mo].

\end{abstract}

\section{Introduction}

In [Ca], Calegari proved that if an atoroidal 3-manifold admits a uniform
1-cochain, then its fundamental group is Gromov-hyperbolic, and it has a
coarse pseudo-Anosov package, which is defined below. These uniform
1-cochains are in some sense a generalization of slitherings, which are
studied by Thurston in [Th]. 

The idea to get the laminations is indeed simple. By [Ca], if M admits a uniform 1-cochain,
M is a Gromov hyperbolic manifold and we have coarse pseudo-Anosov package, 
so there is a coarse lamination in universal cover of the manifold, $\tilde{M}$. 
Using the asymptotic circles of this coarse lamination, we get
a group-invariant family of circles, and using the minimal surface lemmas of
Gabai in [Ga], we can span these circles with laminations by least area
planes. Here we need least area planes to get $\pi_1$ equivariance in universal
cover. Then all we need to show is that this union of laminations in
universal cover can be modified to get a lamination in downstairs, in the
original manifold.\\

\subsection{Definitions:}
The following definitions are from [Ca].

\begin{defn}
{\it Uniform 1-cochain} on a 3-manifold M is a function  
$s:\pi_1(M){\rightarrow}\Bbb R$ satisfying 

\begin{itemize}

\item
$s(\alpha\beta)=s(\beta\alpha)$ for all $\alpha, \beta \in \pi_1(M)$
\item
$s(\alpha^n)=ns(\alpha)$ for any $\alpha\in\pi_1(M)$ and $n \in Z$
\item
$|(\delta s)(\alpha,\beta)|=|s(\alpha)+s(\beta)-s(\alpha\beta)| \leq C_M$, where $C_M$ is a uniform 
constant only depends on M.
\item
For some $t$ the set \[L_t=\{\alpha\in\pi_1(M) \mid |s(\alpha)| \leq t \}\] is coarsely 
connected and coarsely simply connected as a metric space, with the metric inherited as a subspace 
of Cayley($\pi_1(M)$) with some word metric. 

\end{itemize}

\end{defn}

Here, coarsely connected intuitively means that
when you realize $\pi_1(M)$ as a subset of universal cover of M, $\tilde{M}$ (like orbit of a point
under deck transformations), it has an $\epsilon$ neighborhood which is connected, 
and similarly coarsely simply connected means that it has an $\epsilon$ neighborhood 
which is simply connected in $\tilde{M}$.\\

\begin{defn}
{\it A coarse pseudo-Anosov package for M} is the following structure:

\begin{enumerate}

\item
A pair of very full geodesic laminations $\lambda^\pm$ of ${\Bbb H}^2$ which are transverse to each other 
and bind ${\Bbb H}^2$ with transverse measures $\mu^\pm$ without atoms.
\item
An automorphism $Z:{\Bbb H}^2{\rightarrow}{\Bbb H}^2$ which preserves $\lambda^\pm$ and multiplies the measures by k, 
and 1/k respectively.
\item
A uniform quasi-isometry $i:\tilde{M}{\rightarrow}{\Bbb H}^2{\times}\Bbb R$ with the following metric: each level set 
${\Bbb H}^2{\times}n$ is isometric to ${\Bbb H}^2$, and is glued to ${\Bbb H}^2\times(n+1)$ by the mapping cylinder of $Z$ 
whose fibers are normalized to have length 1.
\item
A constant $K$ such that for any $\alpha\in\pi_1(M)$, any $t$, and any $p,q{\in}i^{-1}({\Bbb H}^2{\times}t)$,
$i(\alpha(p))$ and $i(\alpha(q))$ lie on leaves ${\Bbb H}^2{\times}s_1$ and ${\Bbb H}^2{\times}s_2$ where 
$|s_1-s_2| \leq K$.

\end{enumerate}
\end{defn}

This definition might seem awkward at the beginning but, one can think this as a coarse generalization of 
the following structure. Let M be a hyperbolic manifold fibering over $S^1$ with fiber a surface of 
genus greater than 1, $\Sigma$, and the monodromy is pseudo-Anosov map, $\psi$. Then in universal cover, we get a 
${\Bbb H}^2{\times}\Bbb R$ picture as ${\Bbb H}^2$ universal cover of the fiber, $\tilde{\Sigma}$, and 
$\Bbb R$ as universal cover of $S^1$ direction. Now, here we have a pair of lamination $\lambda^\pm$ 
of $\Sigma$ preserved by pseudo-Anosov map, $\psi$. This example fits above definition in the following way:
$\tilde{\lambda^\pm} \subset \hyp^2$ is the very full laminations of $\hyp^2$ in the definition, and $\tilde{\psi}$ is 
the map Z in the definition, and by [CT] there is a quasi-isometry between $\tilde{M}=\hyp^3$ 
and $\hyp^2\times \BR$.

We will call a pseudo-Anosov package {\it transversely orientable} if the lamination $\lambda^\pm$ is transversely 
orientable, and this orientation comes from the $\pi_1(M)$ action on $S^1_\infty({\Bbb H}^2)$. In other words,
$\lambda^\pm$ is transversely orientable lamination and $\pi_1(M)$ action respects this transverse orientation.
Transversely orientable uniform 1-cochain is a uniform 1-cochain which induces transversely orientable 
pseudo-Anosov package.\

\textbf{Notation:}
From now on, 
$\lambda$ will represent a lamination of circle, $S^1_\infty ({\Bbb H}^2)$, 
$\Lambda$ will represent lamination of 3-manifold, 
$\{C\}$ will represent a family of circles in $\Si(\tilde{M})$.
Moreover, if $(x,x')\in S^1 \times S^1$ is an element of lamination $\lambda$, 
$l_x \in \lambda$ will represent corresponding geodesic in ${\Bbb H}^2$ with endpoints
$x,x' \in S^1_\infty ({\Bbb H}^2)$. Similarly, $C_x \in \{C\}$ corresponding circle in $\Si(\tilde{M})$.\\

\subsection{Main Results:}
Our main result is:\\

\textbf{Theorem A:}
Let M be an atoroidal 3-manifold, admitting transversely orientable uniform 1-cochain. Then
there is an induced pair of transverse genuine laminations on M and when you collapse these laminations,
you get a topological pseudo-Anosov flow.\\

\textbf{Outline of the Proof:}

There are 4 main steps:

\begin{enumerate}

\item 
For any leaf $l^+_x\in\lambda^+$ and $l^-_y\in\lambda^-$, we will
assign circles $C^+_x$ and $C^-_y$ in $S^2_{\infty }(\tilde{M})$ such that
the family of circles $\{C^+_x\}$ and $\{C^-_y\}$ are $\pi_1(M)$ invariant
on $S^2_{\infty }$, (i.e. for any $\alpha\in\pi_1(M)$, $\alpha(C^+_x)=C^+_{%
\alpha(x)}$ ).

\item 
We will span this family of circles at infinity, $\{C_{x}\}$ by
laminations of least area planes, $\{\sigma _{x}\}$, such that $\partial
_{\infty }(\sigma _{x})=C_{x}\subset S_{\infty }^{2}$

\item 
We will show that this family of laminations, $\{\sigma_x\}$, are
pairwise disjoint and $\pi_1(M)$ invariant (This is the only step which we use the additional hypothesis of
transverse orientability). Moreover, they induce a pair of genuine laminations $\Lambda^\pm$ on M.

\item 
Using this pair of transverse genuine lamination, we can get a pair of transverse branched surfaces. Then we show 
that this branched surfaces are indeed dynamic pair of branched surfaces which is defined in [Mo]. By [Mo], 
this pair induces a topological pseudo-Anosov flow.

\end{enumerate}

When proving this main theorem, we got very nice by-product. In Step 2 we proved:\\

\textbf{Theorem B:} Let M be a Gromov hyperbolic 3-manifold. Let $\alpha$ be a simple circle in 
$\Si(\tilde{M})$. Then there is a lamination by least area planes spanning this circle $\alpha$ at infinity.

%Moreover, In Step 4, we got:
%
%\textbf{Theorem C:} If $\Lambda^\pm$ are a pair of transverse very full laminations. Then the pair of 
%branched surfaces carrying them is indeed a "dynamic pair of branched surfaces" defined by Mosher in [Mo].
%This means when you collapse these lamination, you get a topological pseudo-Anosov flows.

\subsection{Acknowledgements:} 
I am very grateful to my advisor David Gabai for his patience, guidance, and 
very helpful and inspiring comments. I would like to thank Danny Calegari for his great comments, and 
very helpful explanations. I also thank Sergio Fenley, Tobias Colding and M. Burak Erdogan for very useful conversations.

\section{Preliminaries}

We will give a very rough sketch of some results of Calegari's article [Ca], which is
very crucial for this article. 

Let M be an atoroidal closed 3-manifold, and $s:\pi_1(M)\rightarrow \Bbb R$ be a uniform 1-cochain. 
Let $\tau$ be a "nice" 1-vertex triangulation of M. Consider the lift of $\tau$, 
$\tilde{\tau}\subset \tilde{M}$. Fix a vertex $x_0 \in \tilde{\tau^0}$. Then we can map 
$\pi_1(M)$ to $\tilde{\tau}\subset \tilde{M}$, such that $\alpha \rightarrow 
\alpha(x_0) \in \tilde{\tau^0}$, where $x_0 \in \tilde{\tau^0}$ fixed. 
This is a realization of $\pi_1(M)$ in $\tilde{M}$, as $\pi_1(M) \leftrightarrow \tilde{\tau^0}$ 
with $\alpha \leftrightarrow \alpha(x_0)$. Since $s:\pi_1(M)\rightarrow \Bbb R$, 
we can think of s as a function from a discrete subset
of $\tilde{M}$ to $\Bbb R$. Then extend this function continuously to whole $\tilde{M}$ in a controlled way,
say $S:\tilde{M} \rightarrow \Bbb R$. Now, s is uniform means that, there exist an interval $I\subset \Bbb R$
such that $S^{-1}(I)$ has a k-neighborhood, $N_k(S^{-1}(I))$, which is connected and simply connected.
This is very essential condition as it is used to show that the level sets $\Sigma_t=S^{-1}(t)$ are 
quasi-isometric to $\hyp^2$. 

On the other hand, since $\Sigma_t$ is quasi-isometric to $\hyp^2$, we can talk about the boundary at
infinity of $\Sigma_t$, $\partial_\infty (\Sigma_t) \sim S^1_\infty(\hyp^2)$. The elements  
$x \in \partial_\infty (\Sigma_t)$ are rays, $r_x$, going to infinity. Moreover, he proved that the Hausdorff
distance between any 2 level sets, $d_H(\Sigma_t , \Sigma_{t'})$, is always bounded, and this means 
there is a universal circle $S^1_{univ}$ corresponding to $\partial_\infty (\Sigma_t)$ for any t. In addition,
for any element $\alpha \in \pi_1(M)$, $d_H(\alpha(\Sigma_t), \Sigma_{t+s(\alpha)})$ is bounded by a 
uniform constant.This enables us to define a $\pi_1(M)$ action on $S^1_{univ}$. Let $\alpha\in \pi_1(M)$,
and $x \in S^1_{univ}$, then by using identification $S^1_{univ} \sim \partial_\infty (\Sigma_t)$,
$\alpha(r_x) \in \alpha(\Sigma_t) \sim \Sigma_{t+s(\alpha)} \sim \Sigma_t$ then
$\alpha(r_x) \sim r_y \subset \Sigma_t$, which shows that
$\alpha(x)=y \in S^1_{univ}$ is well-defined.
Then by showing some properties of this canonical action, Calegari got a pair of transverse very full 
measured laminations, $\lambda^\pm$, on $S^1_{univ}$ (which can be thought as geodesic laminations on $\hyp^2$).
Moreover,these measured laminations with a function $Z:\hyp^2 \rightarrow \hyp^2$, 
which preserves $\lambda^\pm$, and expands $\lambda^+$ and contracts $\lambda^-$ gives us a very nice
quasi-metric on $\hyp^2\times \BR$, giving us a quasi-isometric picture of
$\tilde{M}= \bigcup_{t\in \BR} \Sigma_t$ as $\hyp^2\times \BR$. 
By using Bestvina and Feighn's result, he proved that $\pi_1(M)$ is Gromov-hyperbolic. 

Now, we will list some results from [Ca], which we are going to use later:

\begin{itemize}

\item
For any t, t', $d_H(\Sigma_t , \Sigma_{t'}) \leq C_{t,t'}$.

\item
There is a uniform constant C such that for any element $\alpha \in \pi_1(M)$,\\
$d_H(\alpha(\Sigma_t), \Sigma_{t+s(\alpha)}) \leq C$.

\item
$\Sigma_t$ is quasi-isometric to $\hyp^2$

\item
$\pi_1(M)$ acts on $S^1_{univ}$ as described above.

\item
$\pi_1(M)$ action on $S^1_{univ}$, preserves a pair of transverse very full 
measured laminations, $\lambda^\pm$, on $S^1_{univ}$ 

\item 
$\tilde{M}= \bigcup_{t\in \BR} \Sigma_t$ is quasi-isometric to $\hyp^2\times \BR$ with the metric
$ds^2=k^{2t}dx^2 +k^{-2t}dy^2+(logkdt)^2$, where $dx$ represents transverse measure of $\lambda^+$,
$dy$ represents transverse measure of $\lambda^-$, and t is the variable in $\BR$ direction.

\end{itemize}

\subsection{Uniform 1-cochains:}\

3-manifolds admitting uniform 1-cochain are generalizations of 3-manifolds fibering over $S^1$
and 3-manifolds slithering around $S^1$. For example, if M is a 3-manifold fibering over $S^1$,
then let's say $F\rightarrow M \rightarrow S^1$ is the fibration. This induces a map on $\pi_1$
level $s: \pi_1(M) \rightarrow \pi_1(S^1)=\BZ \subset \BR$. This defines a uniform 1-cochain 
except some trivial cases, since the universal cover of the surface F is a plane, and obviously 
coarsely simply connected.

3-manifolds slithering around $S^1$ are generalizations of 3-manifolds fibering over $S^1$. A 3-manifold 
M slithers around $S^1$ if universal cover $\tilde{M}$ fibers over $S^1$ and deck transformations 
respects this fibering, i.e. maps fibers to fibers. If M slithers around $S^1$,
we can induce a uniform 1 cochain for M. Fix a point
$x_0 \in \tilde{M}$, and realize $\pi_1(M)$ in $\tilde{M}$ as the orbit of $x_0$, i.e. 
$\alpha \sim \alpha(x_0)\in \tilde{M}$. Now, if we lift the fibering map 
$f:\tilde{M}\rightarrow S^1$ to $F:\tilde{M}\rightarrow \BR=\tilde{S^1}$, and if we restrict F to 
$\{\pi_1(M)x_0\}$, we get a map $s:\pi_1(M)\rightarrow\BR$. This map does not satisfy the first 2 conditions but 
it satisfies the 3. condition, and using this we can slightly modify our s to satisfy first 2 condition, too.
Define $\overline{s}:=\lim_{n \to \infty} \frac{s(\alpha^n)}{n}$.
Then by using the 3. property, it is easy to check that $\overline{s}$ satisfies the first 3 condition and
Since we slightly modify original s induced from fibering map $S:\tilde{M}\rightarrow\BR$, which has simply
connected fibers, $\overline{s}$ is also uniform 1-cochain on M.

On the other hand, the advantage of the uniform 1-cochains is that they seem very abundant.
If $\pi_1(M)$ is infinite, then $H^1(M)\neq 0$ or geometrization conjecture implies that
second bounded cohomology of $\pi_1(M)$ is nonzero, $H^2_b(\pi_1(M), \BR) \neq 0$, as Gersten proved that
the second bounded cohomology of negatively curved groups are infinite dimensional, see [Ge]. 
This implies that we have lots of bounded 1-cochains satisfying first 3 conditions of uniform 1-cochains.
It might be possible to find some bounded 1-cochains satisfying the topological condition for any manifold
of this kind. Moreover, in [Th], Thurston says that any hyperbolic manifold might be a slithering
around $S^1$ and uniform 1- cochains are coarse generalizations of slitherings.
Because of these reasons, it was believed that they might be all-inclusive class for the hyperbolic part of the geometrization conjecture.
We now know that there are hyperbolic manifolds which are not slitherings. This is because slitherings induce taut foliations, and
there are many hyperbolic manifolds without taut foliations by [RSS]. With this result, we saw that
slitherings are not general enough, so the natural question arised "What about uniform 1-cochains? Are they general enough for weak hyperbolization?".
But the answer was again "No". Last year, Calegari and Dunfield proved that there are also hyperbolic manifolds without uniform 1-cochain, 
by showing Weeks manifold cannot admit uniform 1-cochain, [CD].

\section{Assigning Circles at Infinity}

Now, we will use the following construction of Calegari in [Ca] induced by the
given uniform 1-cochain on M. We will start with a 'nice'
triangulation with one vertex on M, $\tau$. When we lift it to universal cover $\tilde{M}$,
and if we fix a vertex $x_0$ in $\tilde{\tau}$, we can assign each vertex to an
element of $\pi_1(M)$, $\alpha \leftrightarrow \alpha(x_0)$, and 
we get a function from a discrete subset of $\tilde{M}$ to 
$\BR$. We can make a controlled extension so that we get a
function $S:\tilde{M}{\rightarrow}\BR$ induced from the given uniform
1-cochain s. From now on, we fix the unambigiuous triangulation and
controlled extension for M and s. Let $\Sigma_t$ be the level sets of the
function $S:\tilde{M}{\rightarrow}\BR$, i.e. $\Sigma_t=S^{-1}(t)$ for
any $t{\in}\BR$.\newline
Let $t_0{\in}\BR$ fixed.
\\
\\

\begin{lem}
$\partial_\infty(\Sigma_{t_0})=S_\infty^2(\tilde{M})$
\end{lem}

\begin{pf}
Now, by [Ca], we know that ${\forall}t,t'{\in}\Bbb R$, ${\exists}C{\in}\Bbb R$ such that 
$d_H(\Sigma_t,\Sigma_{t'}) {\leq} C$.\\
On the other hand, we know also by [Ca], there is a uniform constant ( independent 
of $\alpha\in\pi_1(M)$ ) such that for any $\alpha\in\pi_1(M)$, 
$d_H(\alpha(\Sigma_t),\Sigma_{t+s(\alpha)}) \leq C$.

Then if we consider the action of $\pi_1(M)$ on $\Si(\tilde{M})$, by above we get
$\alpha (\partial_\infty (\Sigma_t) )=\partial_\infty (\Sigma_{t+s(\alpha)})$. We conclude that for
any $t \in \Bbb R$, for any $\alpha \in \pi_1(M)$, 
$\alpha (\partial_\infty (\Sigma_{t_0}))= \partial_\infty (\Sigma_{t_0})$. Then 
$\partial_\infty (\Sigma_{t_0})$ is $\pi_1(M)$-invariant subset of $\Si(\tilde{M})$. \

Now, let $A=\partial_\infty (\Sigma_{t_0})$, and $C(A)=\bigcup_{x,y \in A}\gamma_{xy}$ where 
$\gamma_{xy}$ represents the geodesic connecting x and y. As A is $\pi_1(M)$-invariant, so is $C(A)$. 
Let $x_0 \in C(A)$, and $B=\{\alpha(x_0)| \alpha \in \pi_1(M)\}$. By invariance of $C(A)$, 
$B \subseteq C(A)$. This implies $\partial_\infty(B) \subseteq \partial_\infty(C(A))$. As M is a 
closed manifold, $\partial_\infty(B)=\Si(\tilde{M})$. The result follows.
\end{pf}

Now, we will recall some notions and results of Cannon-Thurston in the paper [CT]. M is a
3-manifold which is fibering over a circle with fiber a closed surface of
genus 2, S. The monodromy is a pseudo-Anosov map, and so M is hyperbolic
3-manifold.

We are going to make an analogy between this example and our situation.
Consider $S \rightarrow M \rightarrow S^1$ inducing the homomorphism
$s:\pi_1(M) \rightarrow \pi_1(S^1)=Z \subset \BR$. Obviously this is a
uniform 1-cochain for M. So this is a special case of our situation.

We want to analogously extend the following results of [CT]. In the analogy,
we will replace the inclusion of $\hyp^2$ into ${\hyp}^3$, with its coarse correspondent
the inclusion of $\Sigma_{t_0}$ into $\tilde{M}$, and use the result of
[Ca], $\Sigma_{t_0}$ is quasi-isometric to $\hyp^2$.

\begin{itemize}

\item
\begin{equation*}
\begin{array}{ccc}
\BB^2 & \overset{\hat{i}}{\longrightarrow} & \BB^3 \\
\uparrow &  & \uparrow \\
\hyp^2 & \overset{i}{\longrightarrow} & \hyp^3 \\
\downarrow &  & \downarrow \\
S & \hookrightarrow & M%
\end{array}%
\end{equation*}

then i extends continuously $\hat{i}:B^2 \rightarrow B^3$ such that $\hat{i}%
(\partial B^2)=\partial B^3$, or in other words, $\hat{i}(S^1_\infty(\mathbb{%
H}^2))=S^2_{\infty }({\hyp}^3)$, which is a group invariant peano
curve.

\item the diagram
\begin{equation*}
\begin{array}{ccc}
S^1_\infty & \overset{\hat{i}}{\longrightarrow} & S^2_{\infty } \\
p \searrow &  & \nearrow q \\
& S^2 &
\end{array}%
\end{equation*}

commutes, where p is collapsing map of the laminations and q is a
homeomorphism.
\end{itemize}

We are going to prove the above 2 property by following similar techniques of [CT].

\begin{lem}
The inclusion map $i:\Sigma_{t_0} \hookrightarrow \tilde{M}$ extend continuously to
$\hat{i}:S^1_\infty(\Sigma_{t_0}) \rightarrow \Si(\tilde{M})$. Moreover, $\hat{i}$ is $\pi_1(M)$
equivariant.
\end{lem}

\begin{pf}
There are 6 steps:

\begin{enumerate}

\item
$\hat{i}:S^1_\infty(\Sigma_{t_0}) \rightarrow \Si(\tilde{M})$ is $\pi_1(M)$ invariant.\\
The action of $\pi_1(M)$ on $\Si(\tilde{M})$ is defined such that for any point $x$ in 
$\Si(\tilde{M})$, take a ray $r_x$ in $\tilde{M}$ converging to x. Then define 
$\alpha(x)$ as the limit of the ray $\alpha(r_x)$ in $\tilde{M}$.\\

\begin{pf}
Now, since $\partial_\infty(\Sigma_{t_0})=\Si(\tilde{M})$ then for any x in 
$\Si(\tilde{M})$, we can assume $r_x \subset \Sigma_{t_0}$. By the fact that 
$d_H(\Sigma_{t_0+s(\alpha)},\alpha(\Sigma_{t_0})) \leq C$, there exist a ray $s$ in 
$\Sigma_{t_0+s(\alpha)}$ such that $s$ is quasi-isometric to $\alpha(r_x)$. Then by the 
identification between $S^1_\infty(\Sigma_t)$ and  $S^1_\infty(\Sigma_{t'})$ and by the definition 
of action of $\pi_1(M)$ on $S^1_\infty(\Sigma_{t_0})$, this implies the diagram commutes
\end{pf}

\item
For any $x \in S^1_\infty({\Bbb H}^2)$ has arbitrarily small neighborhoods in 
$B^2={\Bbb H}^2 \bigcup S^1_\infty({\Bbb H}^2)$ bounded by closure in $B^2$ of a single leaf of $\{\lambda^+\}$
or $\{\lambda^-\}$.\\

\begin{pf}
By Theorem 6.14 in [Ca], $\{\lambda^\pm\}$ is binding laminations for ${\Bbb H}^2$. Then the result follows 
from Theorem 10.2 in [CT].
\end{pf}

\item
Consider the metric $g$ on $\tilde{M}$ and the $\pi_1(M)$-invariant pseudo-metric 
$ds^2=k^{2t}dx^2+k^{-2t}dy^2+{(logkdt)}^2$ on ${\Bbb H}^2 \times \Bbb R$. then $\varphi_*(ds)$ and g are 
quasi-comparable, where $\varphi:\tilde{M} \rightarrow {\Bbb H}^2 \times \Bbb R$ is the quasi-isometry in the 
coarse pseudo-Anosov package defined in [Ca].\\

\begin{pf}
First, clearly the metric defined in coarse pseudo-Anosov package defined in [Ca], for ${\Bbb H}^2 \times \Bbb R$
is quasi-isometric to the metric $ds^2$, by definition. Now, by theorem 12.1 in [CT] we know, the
metric on ${\Bbb H}^2$ is quasi-comparable to the metric induced by the laminations $\{\lambda^\pm\}$. So, 
$(\tilde{M},g)$ is quasi-comparable to $({\Bbb H}^2 \times \BR, ds^2)$.
\end{pf}

\item
If $l$ is a leaf of $\{\lambda^+\}$ or $\{\lambda^-\}$ in ${\Bbb H}^2$, then $l \times \Bbb R$ is totally 
geodesic in $({\Bbb H}^2 \times \BR, ds^2)$.\\

\begin{pf}
WLOG assume $l$ in $\{\lambda^+\}$. Define 
$\rho:{\Bbb H}^2 \times \BR \rightarrow l \times \Bbb R$ as a product map, $\rho=(f,id)$. Here, 
$f:{\Bbb H}^2{\rightarrow}l$ maps any $l'$ in $\{{\lambda^-}\}$ to $l\cap l'$ (if nonempty), and 
any component $U \subset ({\Bbb H}^2-\lambda^-)$ to $U \cap l$. This retraction is $ds$-reducing as in
Theorem 5.2 in [CT], so $l\times \Bbb R$ is totally geodesic.
\end{pf}

\item
Fix $z\in {\Bbb H}^2$, $\forall\epsilon>0$ $\exists N$ such that if $d_H(z,l) > N$, then the radius of 
$\partial_\infty (l \times \BR) \subset \Si(\tilde{M})$ is less than $\epsilon$.\\

\begin{pf}
The topology is defined as if $a,b \in \Si(\tilde{M})$, and $\gamma_{ab}$ is the geodesic 
connecting $a$ and $b$, then if $\gamma_{ab} \bigcap B_{2^k}(z)= \emptyset$, then 
$d(a,b) < \frac{1}{k}$, by [Gr]. Since, $\tilde{M}$ is negatively curved, then there is a uniform 
constant $C_k$ such that for any k-quasi-geodesic $\alpha_{xy}$ between x and y, 
$d_H(\alpha_{xy},\gamma_{xy}) < C_K$where $\gamma_{xy}$ is the geodesic between x any y, and 
$d_H$ represents Hausdorff distance.Since $l \times \Bbb R$ is quasi totally geodesic in $\tilde{M}$, 
for any $r$, choose $N=2^r+C_k$, where k is the uniform quasi-isometry constant, 
then the radius of $\partial_\infty (l \times \BR) \subset \Si(\tilde{M})$ is less than 
$\frac{1}{k}$ in $\Si(\tilde{M})$.
\end{pf}

\item
Proof of the lemma:\\

Le $x\in S^1_\infty(\Sigma_{t_0})=S^1_\infty({\Bbb H}^2)$, then there exist a sequence of subsets
$C_1\supseteq C_2 \supseteq$ .....$\supseteq C_n\supseteq$ ..... in $B^2={\Bbb H}^2\bigcup S^1_\infty({\Bbb H}^2)$
such that $C_n$ is bounded by a leaf $l_n \in \lambda^\pm$ and $x=\bigcap_{n=1}^\infty C_n$, by step 
(2). Let $U_n={\Bbb H}^2\bigcap C_n$. Define $\hat{i}(x)=\bigcap_{n=1}^\infty \overline{i(U_n)}\subset B^3$

Now, we will prove that $\hat{i}$ is single valued. Consider $i(l_n\times \BR)$ seperates $i(U_n)$ from 
a large compact set. Then by Step (5), as $n\rightarrow\infty$, 
$diam(\partial_\infty(l_n\times \BR))\rightarrow 0$. This means 
$\hat{i}:S^1_\infty(\Sigma_{t_0})\rightarrow\Si(\tilde{M})$ is well-defined. Again, by step (5) 
and above argument, 
$\forall\epsilon>0$, $\exists \delta_n>0$ such that $B_{\delta_n}(x)\subset B^2$ is a neighboorhood 
of x with $\hat{i}(B_{\delta_n}(x)) \subset B_{\epsilon}(\hat{i}(x))\subset \Si(\tilde{M})$.
This proves that $\hat{i}$ is continuous.
\end{enumerate}
\end{pf}

Now, we are going to prove the second property:

\begin{lem}
The Gromov boundary $S^1_\infty(\Sigma_{t_0})$ maps $\pi_1(M)$-equivariantly to a sphere $S^2$, 
by quotienting each leaf in $\{\lambda^\pm\}$ to a point. The quotient sphere is $\pi_1(M)$-equivariantly
equivalent to $\Si(\tilde{M})$.
\end{lem}

\begin{pf}
Again, we will use the method of [CT]. There are 4 main steps. Consider $B^2 \times I$ as 
compactification of ${\Bbb H}^2 \times \Bbb R$

\begin{enumerate}

\item
Extend $\hat{i}:(\partial B^2)=\partial (B^2 \times 0) \rightarrow \partial B^3$ to a map
$\varphi:\partial (B^2 \times I) \rightarrow \partial B^3$.

\item
Define a cellular decomposition G of the 2-sphere $\partial (B^2 \times I)$ 
by using the leaves of the two singular foliations (induced by $\{\lambda^\pm \times \BR \}$ after
collapsing complementary regions), say $F^+$ and $F^-$.

\item
Show that $\varphi$ factors through

\[\begin{array}{ccc}
\partial (B^2 \times I) & \overset{\varphi}{\longrightarrow} & \partial B^3\\
p \searrow &  & \nearrow q\\
 & \partial(B^2 \times I) /G & 
 \end{array}\]
 
 where $\partial(B^2 \times I) /G \simeq S^2$, and G is the decompostion of $\partial(B^2 \times I)$.
 
 \item
 Show that $q:\partial(B^2 \times I) /G \rightarrow \partial B^3$ is homeomorphism. \\
 
\end{enumerate}

{\it Proofs of the steps:}

\begin{enumerate}

\item
Extending $\varphi$:\

We have $\hat{i}:S^1_\infty\rightarrow \Si$. Consider 
$S^1_\infty= \partial (\BB^2 \times \{0\}) \subset \partial(\BB^2 \times I)$.
Now, let $p\in((\partial \BB^2) \times I)$, and let $r_p$ be any ray such that 
$r_p(t)\rightarrow p$ as $t\rightarrow \infty$. If $p\in \hyp^2 \times \{-\infty,+\infty\}$,
then let $r_p$ be the vertical ray asymptotic to p. Then define 
$\varphi(p)=\overline{Q(r_p)}\cap \Si(\tilde{M})$ where $Q:\hyp^2 \times \BR \rightarrow\tilde{M}$ is
the quasi-isometry.\

Now, by the proof of Lemma 2.1, we know, when $p\in ((\partial \BB^2) \times I)$, $\varphi$ is 
well-defined. If $p\in \hyp^2 \times \{-\infty,+\infty\}$, assume $p\in L$, a leaf of $F^+$, ("foliation"),
then since $L\times I$ is totally geodesic, by Lemma 2.2, and it has induced metric 
$ds^2=(k^{-t}dy^2) + {(logkdt)}^2$ since $dx$ is 0 on $L$. By substitution $T=k^t$, we get 
$ds^2= (dy^2 + dT^2)/T^2$, which is hyperbolic plane in half space model. So the vertical
ray is a geodesic. Then, since $({\hyp}^2 \times \BR , ds)$ quasicomparable to $(\tilde{M}, g )$, 
$\overline{Q(r_p)}\cap\Si$ has a single point , as $Q(r_p)$ is quasigeodesic. So, 
%^
$\varphi: \partial (\BB^2 \times I)\rightarrow\Si(\tilde{M}$ is well-defined.\\

\item
Cellular decomposition: \

The cellular decomposition of $\partial(\BB^2 \times I)$ is same with the one in Section 15 of [CT].
There are 3 kinds of element in decomposition:

\begin{itemize}
\item
(First type)
$L\in F^+$,then $g_1=(L\times +\infty)\cup ((\partial L)\times I) \in  G$

\item
(Second type)
$L\in F^-$,then $g_2=(L\times -\infty)\cup ((\partial L)\times I) \in  G$

\item
(Third type)
$p \in (S^1-\lambda^pm)$, then $g_3=p \times I \in G$

\end{itemize}

This decomposition is cellular, as it is proved in section 14 in [CT].

\item
Factoring $\varphi$: \

We show in (1) that $\varphi$ well-defined. 
Now, we want to show that $\varphi$ factors through the decomposition space projection. In other words, 
if G is the cellular decomposition and $g\in G$, then  for any $p,q\in g$, $\varphi(p)=\varphi(q)$.\

if $g$ is the third type, then by the proof of the Lemma 2.2, the result follows.\

if $g$ is the first type, say $g_1=(L\times +\infty)\cup ((\partial L)\times I)$.
Now, consider $L \times I$with the induced metric
$ds^2=(k^{-t}dy^2) + {(logkdt)}^2$, since $dx$ is 0 on $L$.
By substitution $T=k^t$, we get 
$ds^2= (dy^2 + dT^2)/T^2$, which is hyperbolic plane in half space model.
Then consider the geodesics in this space which is in the complement of 
a very large circle, perpendicular to the boundary, say $\gamma_t$ is a geodesic 
which lies in the complement of a radius-$t$ circle. Since $\gamma_t$ is geodesic
in $L \times I$ which is totally geodesic in $\hyp^2 \times \BR$, then 
$\gamma_t$ is also geodesic in $\hyp^2 \times \BR$. This space is quasi-comparable with $\tilde{M}$. 
Hence, $Q(\gamma_t)$ is a quasi-geodesic in $\tilde{M}$, and as $t \rightarrow \infty$, $\gamma_t$ miss 
larger compact sets, then by the definition of the topology in $\Si(\tilde{M})$, the
endpoints of $Q(\gamma_t)$ will converge to a point in $\Si(\tilde{M})$.
This proves that for any $p,q \in g_1$, $\varphi(p)=\varphi(q)$.\

Similar proof works for the second type, too.\\

\item
q is homeomorphism: \

By Theorem 14.1 in [CT], $\partial (\BB^2 \times I) / G \simeq S^2$. Now,

\[\begin{array}{ccc}
\partial (B^2 \times I) & \overset{\varphi}{\longrightarrow}  & \partial B^3\\
p \searrow &  & \nearrow q\\
 & \partial(B^2 \times I) /G & 
 \end{array}\]
 
 By (3), $q$ is well-defined. By Lemma 2.2, $q$ is onto, as $\varphi$ is onto.
 So, we need to show that $q$ is continuous, and injective.\\
 In order to show that $q$ is continuous, it suffices to show $\varphi$ is continuous.
 Since every element in $G$ intersects $\partial \BB^2 \times \{0\}$, then 
 $(\partial \BB^2 \times \{0\}) / G' \simeq \partial (\BB^2 \times I) / G$,
 where $G'$ is the decomposition on $S^1= (\partial \BB^2 \times \{0\})$ induced by $G$.
 So, consider the following commutative diagram:
 
 \[\begin{array}{ccccc}
\partial (\BB^2 \times I) / G & \overset{p_2}{\longleftarrow} & \partial (\BB^2 \times I) & \overset{\varphi}{\longrightarrow}  & \partial B^3\\
\| p_3 &  & \uparrow & \nearrow \hat{i} &  \\
(\partial \BB^2 \times \{0\}) / G' & \overset{p_1}{\longleftarrow} & (\partial \BB^2 \times \{0\}) & & 
\end{array}\]
 
Now, we know $\hat{i}$ is continuous by previous parts. So, for any open set $U \subset \partial \BB^3$,
$\hat{i}^{-1}(U)$ is open in $\partial \BB^2 \times \{0\}$, and since $p_1$ is decomposition space projection
$p_1(\hat{i}^{-1}(U))$ is open in $(\partial \BB^2 \times \{0\}) / G'$. By the homeomorphism,
$p_3(p_1(\hat{i}^{-1}(U))$ is open in $(\BB^2 \times I) / G$ and again since $p_2$
is decomposition space projection, $p_2^{-1}(p_3(p_1(\hat{i}^{-1}(U)))$ is open in $(\BB^2 \times I)$.
Since $\varphi$ factors through $G$, 
$\varphi^{-1}(U)=p_2^{-1}(p_3(p_1(\hat{i}^{-1}(U)))$.
This implies $\varphi^{-1}(U)$ is open, and $\varphi$ is continuous.\

\begin{figure}
\mbox{\vbox{\epsfbox{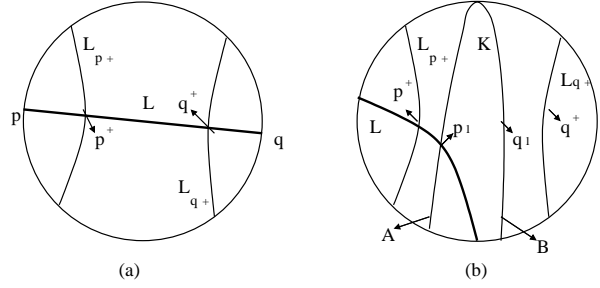}}}
\caption{\label{fig:figure01}
{$p^+$ and $q^+$ are in same leaf $L \in F^+ \times \infty$.}}
\end{figure}

Now, if we show $q$ is injective, then (4) and hence the lemma will be proven.
Clearly, this is equivalent to show, if for any $p,q\in \partial(\BB^2 \times I)$, $\varphi(p)=\varphi(q)$,
then there exist $g \in G$ such that $p,q \in g $.\

Again, we will follow the proof in [CT]. Since $\varphi$ factors through $G$,
we can assume $p,q\in (\partial \BB^2) \times I$.\\

\textbf{Claim 1:}
$\exists p',q' \in \hyp^2 \times \{-\infty, +\infty\}$ arbitrarily close to $p\times I$ and 
$q\times I$ such that $\varphi(p) = \varphi(p') = \varphi(p) = \varphi(q) = \varphi(q')$\\

\textbf{Claim 2:}
$p', q'$ lie in the same element in $G$. \\

Assuming these two claims, we can prove injectiveness as follows. By taking limits,
$p'\rightarrow p$ and $q' \rightarrow q$, we see that $p$ and $q$ are in same element of $G$. 
The result follows. Hence, proving these two claims will be enough.\\

{\it Proof of Claim 1:}
Let $L\times I$ separates $p$ from $q$. Then $L\times I$ separates terminal rays of $r_p$ and $r_q$.
But since $\varphi(p)= \varphi(q)$ then $\varphi(p) \in \partial_\infty(L\times I)$. 
So,we can take $p'\in \partial_\infty(L\times I) \cap (\hyp^2\times \{-\infty, +\infty\})$
such that $\varphi(p)= \varphi(p')$. Since we can choose $L$ close to $p$,
we can assume $p'$ is arbitrarily close to $p$. Similarly for $q$, we can choose $q'$
arbitrarily close to $q'$.\\

{\it Proof of Claim 2:} Let $p^+$ be the projection of $p'$ into $\hyp^2 \times +\infty$,
$p^-$ into $\hyp^2 \times -\infty$. Similarly, define $q^+$, and $q^-$. Let $p'=p^+$.\\

\textbf{Claim :3}
The leaf $L\in F^+$ such that $p^+\in L \times +\infty$ also contains $q^+$, as in Figure [1a].\\

Assuming Claim 3, since $p^+$ and $q^+$ are identified and lie in same $g\in G$, if $q'=q^+$
then we are done. If not, $q' =q^-$ which identified with $p^-$. But we know 
$q'$ and $p'$ are identified by $\varphi$, This means $\varphi(p^+)=\varphi(p^-)$. But we 
know that the vertical geodesic between $p^+$ and $p^-$ is corresponding a quasi-geodesic in $\tilde{M}$,
hence it cannot have only one endpoint at infinity. This establishes Claim 2.\\

{\it Proof of Claim 3:}
$\exists L_{p^+}, L_{q^+} \in F^-\times +\infty$, $p^+ \in L_{p^+}$ and $q^+\in L_{q^+}$.\

First we show $L_{p^+}$ and $L_{q^+}$ are different. Otherwise, $\varphi(p^+)\neq\varphi(q^+)$ (of course 
we are assuming $p^+\neq q^+$) and $\varphi(p^+) \neq \varphi((\partial)\times I)$,  
as $\varphi$ is injective on $F^-\times +\infty$. But $\varphi(q^-) = \varphi((\partial)\times I)$, and
this means $\varphi(p^+)\neq\varphi(q^-)$. This implies $\varphi(p^+)=\varphi(p')\neq\varphi(q')$, which is contradiction.\

Let $L\in F^+ \times +\infty$ such that $p^+ \in L$.
Consider the leaves which separates $p^+$ from $q^+$. Then these leaves 
form an open arc, say $(L_{p^+}, L_{q^+})$ in leaf space of $F^-\times +\infty$. Now, consider the intersection of 
$L$ and the leaves in $(L_{p^+}, L_{q^+})$. If L intersect all of them, and in particular $L_{q^+}$, 
then we are done as $\varphi(L)=\varphi(p^+)=\varphi(q^+)$ and as $\varphi$ is injective on $L_{q^+}$, 
then $L\cap L_{q^+}=\{q^+\}$.
Otherwise, $\exists K \in (L_p, L_q)$ which is the last leaf in 
$(L_p, L_q)$ which L intersects. Then $K$ has maximal subarcs $A$ and $B$ such that $A$ separates 
$p^+$ from $K - A$ and $q^+$, and $B$ separates $q^+$ from $K - B$ and $p^+$, see Figure [1b].
Then as in Claim 1, $\exists p_1\in \overset{\circ}{A} \subset K$ such that 
$\varphi(p_1) = \varphi(p^+)$. Similarly, $\exists q_1\in \overset{\circ}{B} \subset K$ 
such that $\varphi(q_1) = \varphi(q^+)$.
But, since the endpoints of $F^- \times \infty$ is not same with $\partial L$, and $\varphi$ is 
injective on $K$, this implies $p_1=q_1$. But the leaf through any point in $\cirB$ continues into 
a domain of $\hyp^2 \times \infty$ whose closure contains $q_1$ and $q^+$. Then continuation of L 
through $q_1$ intersects further leaves separating $p^+$ and $q^+$. So, K cannot be the last leaf in 
$(L_{p^+}, L_{q^+})$, this is a contradiction. So L intersects $L_{q^+}$ and $\{q^+\}=L\cap L_{q^+}$.
\end{enumerate}
\end{pf}

\begin{figure}
\mbox{\vbox{\epsfbox{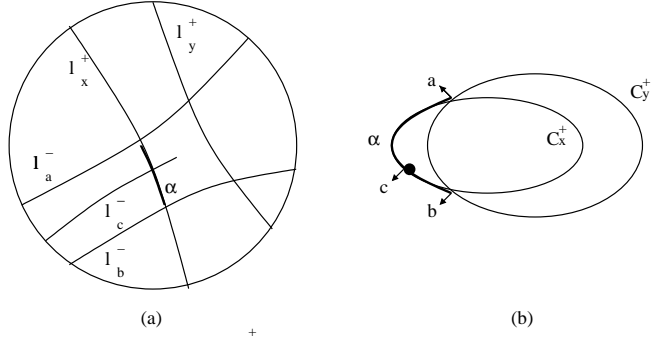}}}
\caption{\label{fig:figure02}
{Intersections of the circles can have at most one component.}}
\end{figure}

\begin{thm}
For any leaf $l^+_x \in \{\lambda^+\}$ and $l^-_y \in \{\lambda^-\}$, there are corresponding 
circles $C^+_x,C^-_y\subset \Si(\tilde{M})$ such that the family of circles $\{C^+_x\}$ and $\{C^-_y\}$
are $\pi_1(M)$ invariant on $\Si(\tilde{M})$, i.e. $\alpha(C^+_x)=C^+_{\alpha(x)}$.
\end{thm}

\begin{pf}
Let $l^+_x \in \{\lambda^+\}$, then consider $l^+_x \times I \subset B^2 \times I$ and 
%^
$\partial_\infty(l^+_x \times I) \subset \partial_\infty(B^2 \times I)$.The collapsing map  
$p:\partial(B^2 \times I) \rightarrow \partial(B^2\times I)/G$ collapses 
$\partial(l^+_x)\times I \bigcup l^+_x \times +\infty$ to a point and maps $l^+_x\times -\infty$
injectively. So, $p(\partial_\infty(l^+_x \times I))$ is a circle in $S^2=\partial(B^2\times I)/G$.
By above, we know that $q:\partial(B^2 \times I) /G \rightarrow \partial B^3$ is homeomorphism.so 
$q(p(\partial_\infty(l^+_x \times I)))$ is a circle in $\Si(\tilde{M})$. Clearly, these circles
are $\pi_1$-invariant by construction.
\end{pf}

\begin{lem}

For any leaf $l^+_x , l^+_y \in \{\lambda^+\}$, the intersection of corresponding 
circles $C^+_x \cap C^+_y\subset \Si(\tilde{M})$ has at most one component, i.e. a point or an interval.

\end{lem}

\begin{pf}
Assume there are more than one component, and choose two points $a,b\in C^+_x \cap C^-_y$ where $a$ and $b$
belongs to different components of intersection.
Consider the proof of previous lemma. We have leaves $l^+_x , l^+_y $ of 
lamination $\{\lambda^+\}$ in ${\Bbb H}^2$, corresponding to the circles. 
Consider the restriction of the map $q\circ p$ to $B^2 \times \{-\infty\}$ and 
the preimages of the points $a$ and $b$ in $B^2 \times \{-\infty\}$. These preimages are going to be two leaves 
$l^-_a , l^-_b \in \{\lambda^-\}$, which are transversely intersecting $l^+_x$ and $l^+_y$, See Figure [2a].
These 4 leaves will define a quadrilateral, $Q$ where each side of belongs to one of them. 
Let $\alpha= Q \cap l^+_x$. $q\circ p(\alpha)$ is a curve in $C^+_x$ which connects $a$ and $b$.
since $a$ and $b$ are in different components of the intersection, There is a point $c$ in $\alpha$ 
whose image is not in $C^+_y$, see Figure [2b]. Then there exist a leaf $l^-_c \in \{\lambda^-\}$ as the preimage of $c$
in $B^2 \times \{-1\}$. Then $l^-_c$ transversely intersect $l^+_x$ but not $l^+_y$. 
But since $l^-_c$ cannot intersect $l^-_a$ and $l^-_b$, then $l^-_c$ cannot go off the quadrilateral $Q$.
This contradicts to fact that leaves are geodesics in ${\hyp}^2 \subset {\BB}^2 \times \{-\infty\}$.
So $C^+_x \cap C^+_y\subset \Si(\tilde{M})$ has at most one component.
\end{pf}

\section{Spanning Circles at Infinity}

We get $\{C^+_x\},\{C^-_y\}\subset S^2_{\infty }(\tilde{M})$, $\pi_1$%
-invariant family of circles at infinity in previous section. Now, we want
to span these circles with laminations by least area planes. If our manifold
were a hyperbolic manifold, then $\tilde{M}\simeq {\hyp}^3$ and the results of Gabai in
[Ga] would give us a positive answer in that situation. But in
our case, the manifold is not hyperbolic, but $\pi_1$-hyperbolic. So, we are
going to extend the results from [Ga], to the case manifold is $\pi_1$%
-hyperbolic. Mainly, we will use the same techniques in Section 3 of [Ga].

\begin{defn} 
If $E\subset B^3= \tilde{M} \bigcup \Si (\tilde{M})$, 
then $C(E)$ denotes the union of geodesics in $\tilde{M}$ connecting points in $E$, i.e. 
$C(E)=\bigcup_{x,y\in E}\gamma_{xy}$ where $\gamma_{xy}$ represents geodesic connecting x and y.
We abuse notation by letting a Riemannian 
metric on $M$ also denote the induced metric on $\tilde{M}$.  An immersed disk with 
boundary $\gamma $ is a {\em least area disc} if it is 
least area among all immersed disks with boundary $\gamma $.  An 
injectively immersed plane is a {\em least area  plane} if each compact subdisk is a least area disk.  

A codimension-$k$ lamination $\sigma $ in the $n$-manifold $Y$ is a 
codimension-$k$ foliated closed 
subset of $Y$, i.e. 
$Y$ is covered by charts of the form $R ^{n-k}\times \BR ^{k}$ and 
$\sigma \mid R ^{n-k}\times \BR ^{k}$ is the product lamination on $R
^{n-k}\times C$,
where $C$  a closed subset
of $R ^{k}$.   Here and later we abuse notation by letting $\sigma $ also 
denote the {\em underlying
space} of its lamination, i.e. the points of $Y$ which  lie in leaves 
of $\sigma $.  Laminations in
this paper  will be codimension-1 in manifolds of dimension 2 or 3.

A {\em complementary region} $J$ is a component of $Y-\sigma $.  Given a 
Riemannian metric on $Y$, $J$ has
an induced path metric, the distance between two points being the 
infimum of lengths of paths in
$J$ connecting them.  A {\em closed complementary region} is the metric 
completion of a complementary
region with the induced path metric.  As a manifold with boundary, a 
closed complementary region is
independent of metric.
\end{defn}

\begin{defn} The sequence $\{S_{i}\}$ of embedded 
surfaces or laminations in a
Riemannian manifold $Y$ converges to the  lamination $\sigma $ if

ia) $\sigma =\{x= $ Lim$_{i\to \infty }x_{i}\mid x_{i}\in S_{i}$ and $\{x_{i}\}$
a  convergent
sequence in $Y\}$;

ib) $\sigma=\{x= $ Lim$_{n_{i}\to \infty }x_{n_{i}}\mid \{n_{i}\}$ an
increasing sequence in ${\mathbb{N}}, x_{n_{i}}\in S_{n_{i}}$ and
$\{x_{n_{i}}\}$ a  convergent
sequence in $Y\}$ $\overset {\text{def}}{=}$ Lim$\{S_{i}\}$. 

ii) Given $x, \{x_{i}\}$ as above, there exist embeddings $f_{i}:D^{2}\to L_{x_{i}}$ which converge in the 
$C^{\infty }$-topology to a smooth embedding $f:D^{2}\to L_{x}$, where 
$x_{i}\in f_{i}(\overset{\circ}{D} {}^{2})$,
$L_{x_{i}}$ is the leaf of $S_{i}$  through $x_{i}$, and $L_{x}$ is the leaf of 
$\sigma $
through $x$, and
$x\in f(\overset{\circ}{D} {}^{2})$. 
\end{defn}

\begin{lem}
If $\{S_i\}$ is a sequence of least area disks in $\tilde{M}$, where $\partial S_i\rightarrow \infty$,
then after passing to a subsequence $\{S_i\}$ converges to a (possibly empty) lamination by least area planes.
\end{lem}

\begin{pf}
There are 5 main steps:

\begin{enumerate}

\item
After passing to a subsequence $Lim\{S_i\}=\{x=\lim_{i \to \infty }x_{i}\mid x_{i}\in S_{i}$ and $\{x_{i}\}$
a  convergent sequence in $\tilde{M}\}$ is closed.\\

\begin{pf}
For each $j$ subdivide $\tilde{M}$ into finite number of closed regions, such that the $j+1$'st subdivision
restricted to B converges to 0, for any compact ball B. In other words, $\tilde{M}=\bigcup_{k=1}^{n_j}B^j_k$
where $B^{j-1}_i=B^j_{i_1} \bigcup$...$\bigcup B^j_{i_r}$, and for compact B $diam(B\bigcap B^j_{n_j}) \rightarrow 0$
as $j\rightarrow \infty$. Now, choose a subsequence of $\{S_i\}$ such that if $i\geq j$ and 
$S_i\bigcap B^j_r\neq\emptyset$, then for any $k>i$, $S_k\bigcap B^j_r \neq \emptyset$. For this subsequence
the limit set $Z=Lim\{S_i\}$ is closed, as for any subsequence in Z, you can use diagonal sequence argument to prove $lim z_i\in Z$.
\end{pf}

\item
Let $\{z_i\}$ be a countable dense subset of $Z$. $\exists\epsilon>0$ such that after passing to a subsequence
of $\{S_j\}$ the following holds. For any $i$, there exists a sequence of embedded disks $D^i_j\subset S_j$ 
which converges to a smoothly embedded least area disk $D_i$ such that $z_i\in D_i$ and 
$\partial D \bigcap B_{\epsilon}(z_i)=\emptyset$.\\

\begin{pf}
Since M is compact we can assume that $\exists\epsilon>0$ such that for any $x\in \tilde{M}$, 
$B_{2\epsilon}(x)$ has strictly convex boundary. Now, fix i, then if 
$D^i_j\subset S_j \bigcap B_{2\epsilon}(z_i)$ is a component,
then $d(z_i, D^i_j) \rightarrow 0$ as $j\rightarrow \infty$. Since $D^i_j$'s are least area, for any j, 
$Area(D^i_j)\leq \frac{1}{2}Area(\partial B_{2\epsilon}(z_i))$. Then by Lemma 3.3 in [HS], 
after passing to a subsequence and 
resricting to $B_{\epsilon}(z_i)$, $D^i_j$'s converge to the desired disk $D_i$. Since this is true for 
each i, the diagonal sequence argument completes the proof.
\end{pf}

\item
There is a lamination $\sigma$ with underlying space Z, such that each $D_i$ is contained in a leaf.
Furthermore $\{S_i\}$ converges to $\sigma$. \\

\begin{pf}
By Step 1, i) of Definition 3.2 holds.  By 
Step 2, for each $i$, $D_{i}\subset Z$.  If 
$x\in \overset\circ{D} _{i}\cap \overset\circ{D} _{j}$, then $D_{i}$ and $D_{j} $ coincide in a 
neighborhood of $x$.  Otherwise being minimal 
surfaces, $D_{i}$ and $D_{j}$ would cross transversely at some point
close to $x$, which would imply that $S_{k}$ was
not  embedded for $k$ sufficiently large, by Lemma 3.6 of [HS].  If $z\in Z$, then the
argument of Step 2 shows that there exists a  convergent sequence
$\{D_{z_{i}}\}\to D_{z}$, where $D_{z_{i}}$ is a subdisk of some $D_{j}$,
$z\in D_{z}$ and $\partial D_{z}\cap B_\epsilon(z)=\emptyset $.   Again
since the $D_{i}$'s  pairwise either locally coincide or are disjoint, 
$D_{z}$ is uniquely determined in an $\epsilon $-neigborhood of $z$. 
Thus $Z=\bigcup _{z\in Z}D_{z}$.   Using the  $D_{z}$'s to define a
 topology on $Z$, it follows that connected components are leaves of a 
lamination  $\sigma $ with
underlying space $Z$.  The uniqueness of $D_{z}$ in $B_\epsilon(z)$ 
implies that near $z$ leaves of
$\sigma $ are graphs of functions over $D_{z}$ and that $\{S_{i}\} $ 
converges to $\sigma $.
\end{pf}

\item
If $g:D\to L$ is an immersion of a disk into a leaf $L$ 
of $\sigma $, then for
all $i$  sufficiently large there exists an immersion $g_{i}:D\to S_{i}$ 
such that $g_{i}\to g $ in the
$C^{\infty }$  topology.\\

\begin{pf}
This is true as $\{S_{i}\}$ converges to $\sigma $.
\end{pf}

\item
Each leaf $L $ of $ \sigma $ is a least area
plane.\\

\begin{pf}
First, we will prove L is a plane.
Let $\tau $ be an essential simple closed curve 
in $L$ and $A\subset L$ a
thin (e.g. $<.5\epsilon $) regular neighborhood of $\tau $.  Let 
$B\subset \tilde{M}$ be a 3-ball transverse to $\bigcup S_{i}$ such that 
$A\subset \overset\circ{B} $.  Let $g:D\to L$ be an isometric immersion of a disc
such that $g(D)=A$ and Area$(D)>$Area$(\partial B)$. (Think of $D$
as being a long thin rectangle.)  By Step 4, for $i$ sufficiently
large, $g$ is closely approximated by an isometric immersion of a
2-disc, i.e. $g_{i}:D_{i}\to S_{i}$ and Area$(D_{i})>$Area$(\partial B)$. 
For $i$ sufficiently large $g_{i}(D_{i})$ is an annulus which closely
approximates $A$.  Otherwise  $g_{i}(D_{i})$ is an embedded disk which
spirals around and closely approximates $A$.     This contradicts
the fact that if $B$ is a ball and $\partial S_{i}\cap B=\emptyset $,
then Area$_{r}(P)\le 1/2$Area$_{r}(\partial B)$, where $P$ is a
component of $S_{i}\cap B$.  Thus for each sufficiently large $i$,
there exists an embedded simple closed curve $\tau _{i}\subset S_{i}$
such that $\{\tau _{i}\}$ converges to $\tau $.  Each $\tau _{i}$ bounds a
disk $E_{i}\subset S_{i}$ of uniformly bounded area.  The
sequence of disks $\{E_{i}\}$ converges to a disk in $L$ bounded by
$\tau $ via  arguments similar to those of the proof of Step 3.  Thus
$L$ is simply connected.  $L$ is not a sphere else for $i$
sufficiently large each $S_{i}$ would be a sphere.\\
Since each embedded  subdisk of $L$ is the limit of least area disks by Step
4,  each embedded subdisk of $L$ is least area and hence $L$ is a
least area plane.
\end{pf}

\end{enumerate}

\end{pf}

\begin{defn}

Let $\alpha $ be an unknotted simple closed 
curve in $\tilde{M}$ with the $r$-metric. Change the $r$-metric of 
$U=\tilde{M}-\overset{\circ}{N} (\alpha )$ by one which coincides with 
$r$ away from a very small neighborhood of $\partial U$ and which gives $U$ a strictly 
convex boundary. It follows by [MSY] that an essential  simple closed 
curve on $\partial N(\alpha )$, also called $\alpha $, 
bounds a properly embedded disk $D\subset U$, least area among all 
immersed disks $E\subset U$ with 
$\partial E\subset \partial U$ and $\partial E$ essential in $\partial U$.   
Call a disk that arises from this construction 
a {\em relatively  least area}  disk in 
$\tilde{M}$

\end{defn}

\begin{lem}
Let $r_{t}$ be a $[0,1]$-parameter family of 
Riemannian metrics on 
$\tilde{M}$ obtained by lifting a $[0,1]$-parameter family on a closed 
 manifold $M$. There exists $e>0$ such
that if $S$ is a relatively least area   disk in $\tilde{M}$ with the 
$r_{i}$-metric, 
then $S\subset N_{\rho }(e,C(\partial S))$
\end{lem}

\begin{pf} 
A short outline: Assume  there is no such e. Then there exists a
sequence of disks $\{D_i\}$ and $D_i\rightarrow \tilde L$ a least area plane
such that $\partial_\infty \tilde L=x$. Moreover, we can choose this $\tilde L$ as $\pi_1$-invariant in 
$\tilde M$. When we project $\tilde L$ to $M$, we see that L is a leaf of an essential lamination
by least area planes. But this implies $M\simeq T^3$ by [Imanishi].

There are 4 main steps:

\begin{enumerate}

\item
There exists an $r$-least area 
plane $\tilde L$ which is a
leaf of a $D^{2}$-limit lamination, and $\partial_\infty \tilde L=x$, where 
$x\in \Si (\tilde{M})$.\\

\begin{pf}
Suppose that for
each $i$, there exists a  relatively $r_{i}$-least area disk $D_{i}^{\prime }$ 
such  that
$D_{i}^{\prime }\not \subset N_i(C(\partial D_{i}^{\prime }))$, where $C(\partial D_{i}^{\prime })$
is the union of geodesics between points in $\partial D_{i}^{\prime }$.
Let $z_{i}\in D_{i}^{\prime }$ be a point 
farthest
from $C(\partial D_{i}^{\prime })$.  A covering transformation of 
$q:\tilde{M}\to M $ 
is an 
isometry in both the $r_{i}$ and $r$ metric.  Therefore by 
replacing each $D_{i}^{\prime }$ by a covering
translate and passing to a subsequence, we can assume that the $z_{i}$ 
converge to  fixed $z_0\in \tilde{M}$. 
By passing to another
subsequence  we  can assume
that Lim$\{C(\partial (D_{i}^{\prime }))\}=w\subset \Si $. Otherwise, it would contain 2 points,
say $x,y\in \Si(\tilde M)$, then $\exists \gamma_{xy}\in\tilde M$.  By using this, 
we can find an upper bound for $d(z_0,C(\partial D_{i}^{\prime }))$. There are sequences $\{x_i\}$ and 
$\{y_i\}$ in $C(\partial D_{i}^{\prime })$, and so there are geodesics $\gamma^i_{xy}$ in 
$C(\partial D_{i}^{\prime })$.
As M is negatively curved, we can get an upper bound for $d(z_0, C(\partial D_{i}^{\prime }))$, 
which is a contradiction. We can cut down the size of the relatively least area disks and pass to a 
subsequence of least area disks $\{D_i\}$. Then by previous lemma, after passing to a subsequence,
we get $D_i\rightarrow \sigma$, where $\sigma$ is the lamination by least area planes. Let $\tilde L$ 
be the leaf containig $z_0$. Replace $D_i$ with $B_i(z_o)\bigcap \tilde L$.
\end{pf}

\item
Let $G_{M}$ denote the group of covering translations of 
$\tilde{M}$ associated to
$M$.   There exists an $r$-least area plane $\tilde Q$ such that for 
each $g\in G_{M}$, either $g(\tilde Q)=\tilde Q$ or $g(\tilde Q)\cap \tilde Q=\emptyset $.  Furthermore  either 
$g(\tilde Q)\cap \tilde L= \emptyset $ or $g(\tilde Q)=\tilde L$.  \\

\begin{pf} There are 2 cases. \\

\textbf{Case 1:} If $w$ is not 
the fixed point of any element of $G_{M}$, then $\tilde L$ is the desired 
$\tilde Q$, otherwise there exists $g\in G_{M}$ such that $g\neq $ id and 
$g(
\tilde L)\cap \tilde L\notin \{\emptyset ,\tilde L\}$.  Since $g(w)\neq w$,
 there exists some 
$i$ such that $g(D_{i})\cap D_{i}\neq \emptyset $
but $g(\partial D_{i})\cap (\partial D_{i})=\emptyset $.  This leads to a
contradiction by the exchange roundoff trick.\\

\textbf{Case 2:} If $w$ is a fixed point of an element of $G_{M}$.\\

We need a lemma for Gromov hyperbolic manifolds, corresponding the fact that
the fundamental group of a closed hyperbolic manifold has no parabolic elements.

\begin{lem}
If M is a closed $\delta$-hyperbolic manifold, every f in $\pi_1(M)$ has 2 fixed point in gromov sphere at infinity.
\end{lem}

\begin{pf}
Assume f has more than 2 fixed points.Let $a, b, c\in\Si$ be fixed points of f. Consider geodesic between
a and b, $\gamma_{ab}$. Since a and b are fixed points of f, $f(\gamma_{ab})=\gamma_{ab}$, this is also
true for $\gamma_{bc}, \gamma_{ca}$. Since there is no fixed point in $\tilde{M}$, f must iterate these 3 geodesics. 
WLOG assume F iterates $\gamma_{ab}$ from a to b, and $\gamma_{bc}$ from b to c. Now, let's take a point 
$x\in\gamma_{ab}$, and another point $y\in\gamma_{bc}$. Now consider geodesic segment between x and y.
Since f is isometry of $\tilde{M}$, the length of [x,y] must be same with the length of $f^n([x,y])$. 
But, since $f^n(x)\rightarrow b$ and $f^n(y)\rightarrow c$ , the length of $f^n([x,y])$ must go to infinity,
so this is a contradiction. This means f  cannot have more than 2 fixed points in $\Si$.\

Now, we will show that f cannot have only one fixed point in $\Si$. This is actually analogous with that closed 
hyperbolic manifolds cannot have parabolic hyperbolic isometries in deck transformations. Assume $a\in\Si$
is the only fixed point of f. Let $b\in\Si$ be an arbitrary point and $c=f(b)$. Consider geodesics 
$\gamma_{ab},\gamma_{ac}$. Let $x$ be an arbitrary point in $\gamma_{ab}$, and $y=f(x)\in\gamma_{ac}$
parametrize geodesics by arclength so that $\gamma_{ab}(0)=x$ and $\gamma_{bc}(0)=y$ with 
$\gamma_{ab}(t)\rightarrow a$ and $\gamma_{ac}(t)\rightarrow a$ as $t\rightarrow\infty$. Then since f is isometry 
$f(\gamma_{ab}(t))=\gamma_{ac}(t)$. But since $\tilde{M}$ $\delta$-hyperbolic, geodesics diverge exponentially the distance between 
$\gamma_{ab}(t)$ and $\gamma_{ac}(t)$ will decrease, that means as $t\rightarrow \infty$ $d(\gamma_{ab}(t),\gamma_{ac}(t))\rightarrow 0$.
But since M is closed there is no cusps, so the length of essential loops is bounded below. This is a contradiction.
\end{pf}

Let $w$ be the fixed point of some primitive 
element $f$ of $G_{M}$. We
find  $Q$ as follows.   Let $A_{f}$ denote the  
axis of $f$.  There does not exist $N>0$ such that 
$\tilde L\subset N_N (A_{f})$.  Otherwise, if $\tilde{L} \in N_{N_0} (A_f)$ then for any
$x\in A_f$ $Area(H_x\cap \tilde{L})\leq \frac{1}{2}Area(\partial B_{N_0})$ where 
$H_x \subset N_{N_0} (A_f)$ cut by a disk in $B_{N_0(x)}$. But, this contradicts to monotonicity formula
(Lemma 2.3. [HS]) as $x\rightarrow w$, the intrinsic radius of the region enclosed by 
$H_x\cap\tilde{L}$ in $\tilde{L}$ goes to infinity whereas the area remains bounded.\

Let $\{y_{i}\}$
be a sequence of points in  $\tilde L$
such that $d(y_{i},A_{f})>i$.  Let $g_{i}\in G_{M}$ be such that  $g_{i}(y_{i})=v_{i}$ 
lies in a
fixed $X$-fundamental domain $V$ in $\tilde{M}$.    By passing to a 
subsequence we
can  assume that $v_{i}\to v\in \tilde{M}$ and $g_{i}(w)\to w'$.  By passing 
to another subsequence we can assume that 
$i\neq j$ implies that $w_{i}\overset{\text{def}}{=}g_{i}(w)\neq g_{j}(w)
\overset{\text{def}}{=}w_{j}$.  Suppose on the 
contrary that for all $i,j, g_{i}(w)=g_{j}(w)$.  
Then $g_{i}(w)=g_{j}(w)\implies g^{-1}_{j}\circ g_{i} (w)=w\implies g^{-1} _{j}\circ g_{i}=f^{n_{i}}\implies g_{i}=g_{j}\circ f^{n_{i}}$.  
Now $g_{i}(y_{i})\subset V \implies y_{i}\in g^{-1} _{i}(V)=f^{-n_{i}}\circ g^{-1} _{j}(V)\implies d(y_{i},A_{f})\le $max$\{d(g^{-1} _{j}(z),A_{f})\mid z\in V\}$.   The finiteness of the latter contradicts the choice of 
$y_{i}$, for $i$ large.

Let $\tilde Q$ be a least area plane passing through $v$, obtained by  applying
Lemma 4.1. to the sequence  $g_{i}(\tilde L)=\tilde L_{i}$, or more precisely to
$\{g_{i}(D_{n_{i}})\}$,  where $\{n_{i}\}$ is a sufficiently fast  growing
sequence. There exists no $h\in G_{M}$ such  that 
%^
$h(\tilde Q)\cap \tilde
Q\notin \{\emptyset ,\tilde Q\}$;  else for sufficiently large $i,j$, 
$h(\tilde L_{j})\cap \tilde L_{i}\neq \emptyset $.  Therefore there exists 
$i,j$ such that $h(\tilde L_{j})\cap \tilde L_{i}\neq \emptyset $ and
$w_{i}\neq h(w_{j})$.   This implies that $g^{-1} _{i}\circ h\circ  g_{j} 
(\tilde L)\cap \tilde L\neq \emptyset $ and $g^{-1} _{i}\circ h\circ
g_{j}(w)\neq w$, which is a contradiction.  A similar argument shows that
$h(\tilde L)\cap \tilde Q\in \{\emptyset ,\tilde Q\}$.
\end{pf}

\item
There exists a least area properly embedded
plane $\tilde P$ with $\partial_\infty (\tilde P)$ is a point in $\Si(\tilde M)$ such that
for each  $g\in G_{M}, g(\tilde P)=\tilde P$ or  $g(\tilde P)\cap \tilde P=\emptyset $.  
If $\pi :\tilde{M}\to M$ is the covering
projection, then  $\pi (\tilde P$) projects to a  leaf $P$ of an
essential  lamination $\kappa $ in $M$.  Finally the leaves of
$\kappa $ lift to least area planes in $\tilde{M}$ and each
leaf of $\kappa $ is dense in $\kappa $.\\

\begin{pf}
Let $\lambda $ be the lamination in $X$ obtained 
by taking the 
closure of 
the injectively immersed surface $Q$ which is the 
projection of $\tilde Q$.   We show that $\lambda $ is essential by 
showing 
that each leaf is incompressible 
and end incompressible [GO].   Each leaf $Q_{\alpha }$ of $\lambda $ lifts 
to a 
surface 
$\tilde Q_{\alpha }$ in $\tilde{M}$ which is a limit 
of translates of subdisks of $\tilde Q$, hence $\tilde Q_{\alpha }$ is a 
leaf of a $D^{2}$-limit
lamination and hence is a least  area plane, 
so 
$Q_{\alpha }$ is incompressible.  An end 
compression of $Q_{\alpha }$ would imply the existence  of a monogon in
$\tilde{M}$ connecting two very  close together subdisks of $\tilde Q$ of
very much larger area,
contradicting the fact that $\tilde Q_{\alpha }$ is least area as in  
Figure [3].

\begin{figure}
\mbox{\vbox{\epsfbox{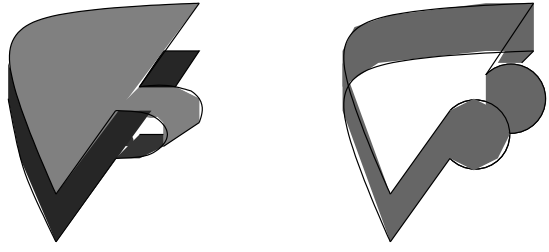}}}
\caption{\label{fig:figure03}
{least area planes are end-incompressible.}}
\end{figure}

  Let 
$\kappa $ be a  nontrivial sublamination of $\lambda $ such 
that each leaf of $\kappa $ is dense in $\kappa $.  

The lift 
$\tilde \kappa $ of $\kappa $ to $\tilde{M}$ is a sublamination of the 
lamination 
which is the closure of all the  $G_{M}$-translates of $\tilde Q$.  Since
$\tilde L $ is either  disjoint from $\tilde \kappa $ or a leaf of $\tilde \kappa $, it follows that $L=\pi (\tilde L)$ is either 
a leaf of
$\kappa $ or disjoint from $\kappa $. By construction $\kappa \subset \overline L$ since $\tilde Q$ is in the closure of
$G_{M}(\tilde L)$.

If $\tilde L$ is a leaf of $\tilde \kappa $, then Step 3 holds with 
$\tilde P=\tilde L$.  In that case since $\tilde L$ is the lift of a 
leaf of an essential lamination, it
follows by [GO] that $\tilde L$ is properly embedded in $\tilde{M}$.\

Now, we will show, that if $L\subset J$, where $J$ is a complementary region of $\kappa$, then 
L can be replaced with a leaf of the foliation, say P, which lies in the boundary of the complementary
region $J$, and P has the desired properties.\\

\textbf{Claim:}
$J=\overset{\circ}{D} {}^{2}\times I$ and  $L$ is homotopic to
$\overset{\circ}{D} \times 1/2$ via a homotopy in
$J$ in which  points of $L$ are moved by  homotopy tracks of uniformly
bounded length. \\

{\it Proof of Claim:}
As in [GO] $J$ is of the
form $A\cup Z$, where each component of interstitial region $A$ is an $I$-bundle over a
noncompact surface, gut region $Z$ is a
connected compact $3$-manifold and $A\cap Z$ is a union of annuli.  Since
$M$ is of finite volume,
by taking $Z$ to be sufficiently big (by reducing the size of $A$) we
can assume that the
$I$-fibres are very short $\rho $-geodesic arcs nearly orthogonal to
$\partial J$.  Since L is least area plane which means it is tight in some sense
(by [S], L has bounded second fundamental form) if the $I$-fibres
are sufficiently short, then they must be transverse to $L$.  Thus we
can assume that $L$ is transverse to the $I$-fibres of $A$.

Assume $A\neq \emptyset$. If $E$ is a vertical annulus in $A$, i.e.  a union of $I$-fibres,
then either $E$ spans a $D^{2}\times I\subset J$ or
$E\cap L=\emptyset $.  Otherwise $E$ lifts to an $I\times \Bbb R $ whose core
$\alpha $ is properly homotopic (by the previous paragraph) to a
curve lying in $\tilde L$, contradicting Step 1, for $\alpha $ has
distinct endpoints in $\Si $.    Since $\kappa \subset \overline L$, it
follows that  some component $A_{1}$ of $A$ and hence  each component
of $A_{1}\cap Z$ nontrivially intersect $L$ and hence $A_{1}=A$ and $J$
is obtained by attaching 2-handles to $A$ along $A\cap Z$.  Since
each vertical annulus in $A$ bounds a $D^{2}\times I$, it follows that
$J=\overset{\circ}{D} {}^{2}\times I$. Since $J$ is simply connected, it lifts to
$\tilde{M}$ and hence $L$ is embedded in $J$ since $\tilde L$ is
embedded in $\tilde{M}$.  Therefore  if $E\subset A$ is a vertical
annulus, then $E\cap L$ is a union of embedded circles.  Each
such circle bounds a disk in $L$ which is isotopic rel boundary to a
horizontal disk in the associated $D^{2}\times I$.  If $P$ is a
component of $\partial J$, then  vertical projection of $L\cap A$ to
$P\cap A$ extends to an immersion of $L$ to $P$.  $P$ being simply
connected implies that this is in fact a diffeomorphism.  Again as
in [GO] each lift of $P$ is properly embedded.

If $A=\emptyset$, derive a contradiction as follows. In this case $\kappa$ is a
closed $\pi_1$-injective surface $S_0$. Consider an incompressible surface
$S_1$ in $X$ split open along $S_0$ which nontrivially intersects $S_0$ and
consider $L\cap S_1$ to argue that the limit set of $\tilde L$ consists of more
than a point.

Since each leaf of $\kappa$ is dense in $\kappa$ the above argument shows that
$\kappa$ has no closed leaves.
\end{pf}

\item

Proof of Lemma.\\

\begin{pf}
Note that $\tilde P$ could have been chosen so that $w \in \Si(\tilde M) $ is
the asymptotic boundary of $\tilde P$, $\partial_onfty(\tilde P)$.
If $B$ is the region in $\BB ^{3}$ bounded by $\tilde P$
such that $B\cap \Si =w$, then $G_{B}=\{g\in G_{M}\mid g(\cirB )\cap \cirB \neq \emptyset \}$
is a subgroup of the stabilizer $G_{w}$ of $w$.  Since
$G_{w}$ is generated by $f$, $G_{B}$ is generated by $f^{n}$ for some
$n\in \BZ $.   First suppose that $G_{B}\neq $ id.  We can assume
that   $f^{n}(B-w)\subset \cirB $. Since $\tilde P$ is proper, each
$z\in \tilde P$  has a neighborhood $W\subset \tilde M$ such that
$W\cap (f^{n}(\tilde P)\cup f^{-n}(\tilde P))=\emptyset $ and hence
$\{g\in G_{M}\mid g(\tilde P) \cap W\neq \emptyset \}=$
id.  This implies that $P$ is isolated, contradicting the fact that
each leaf of $\kappa $ is
dense and $\kappa $ has no closed leaves.  Finally consider the case
$G_{B}=$ id. In this case $\cirB \cap \tilde \kappa =\emptyset $,
otherwise $P$ is dense in $\kappa $ implies that some covering
translate of $\tilde P$ lies in $\cirB $.    Let $I$ be an $I$-fibre
of $A$ and let $\tilde I$ be the lift which intersects $\tilde P$.
Since $P$ is nonisolated,  $\tilde I\subset B$, with one endpoint
$i\in \cirB $.  We obtain the contradiction
$\tilde \kappa \cap \cirB \neq \emptyset $.
\end{pf}

\end{enumerate}

\end{pf}

\begin{thm}
Let
$\tau $ be a
simple closed curve in $\Si $.  Then
there exists a $D^{2}$-limit
lamination $\sigma \subset \tilde{M}$ by least area planes spanning
$\tau $.
Furthermore there exists $e>0$,(independent of $\tau $),  such
that if $\sigma $ is any spanning lamination  by least area 
planes, then $\sigma \subset N_e(C(\tau ))$.

\end{thm}

\begin{pf}
Let $e>0$ be as in Lemma 3.7.  Let $\omega $ be a
properly embedded path in $\BB ^{3}$
connecting points in distinct components of $\Si -\tau $.
We will prove this lemma in 5 steps.

\begin{enumerate}

\item
$N_{5e}(C(\tau)) \simeq \overset{\circ}{D} {}^{2} \times I$\\

\begin{pf}
Let $\gamma_{xy}$ be the geodesic between x and y, where $x,y\in\Si$. Let $D_x:=\bigcup_{t\in\tau} \gamma_{xt}$. 
Then $C(\tau)=\bigcup_{x\in\tau} D_x$. We first prove that $N_{2\delta}(D_x)\simeq \overset\circ{D^2}\times I$. Fix $t_0\in\tau$.
Let $\{t_n\}\subset \tau$ and $t_n\rightarrow t_0$. Let $a_n\in\gamma_{xt_0}$ such that 
$\gamma_{xa_n}=\gamma_{xt_0} \cap N_{2\delta}(\gamma_{xt_n})$. Since $\tilde {M}$ is $\delta$-hyperbolic, the triangles are $\delta$ thin, 
if vertices are in $\tilde{M}$, and 2$\delta $ -thin if the vertices in $\Si(\tilde{M})$. So as $t_n \rightarrow t_0$, $a_n \rightarrow t_0$
That means $N_{2\delta}(D_x)= \cup_{t\in\tau} N_{2\delta} \gamma_{xt}$, so it is homeomorphic to $\overset{\circ}{D^2}\times I$.

Now, consider $C(\tau)=\bigcup_{x\in\tau} D_x$. $\forall x,y\in\tau$, $d_H(D_x,D_y) < 2\delta$ since for any $u\in D_x$, 
$u\in\gamma_{xt}$ then since $\gamma_{xt}\cup\gamma_{xy}\cup\gamma_{yt}$ is a 2$\delta$-thin triangle, so $\gamma_{xt}\subset N_{2\delta}\gamma_{yt}\cup\gamma_{xy}$
this means $u\in\gamma_{xt}\subset N_{2\delta}(D_y)$. This proves $d_H(D_x,D_y) < 2\delta$.That shows $D_x\subset N_{2\delta}(D_y)$, so 
$ N_{2\delta}(C(\tau))=\bigcup_{x\in\tau} N_{2\delta}(D_x)$ is homeomorphic to $\overset{\circ}{D^2} \times I$.
If $ e < 2\delta$, replace e such that $e > 2\delta$, then result follows.
\end{pf}

\item
There exists a sequence of  relatively least area
disks $\{E_{i}\}$ such that for  each $i, E_{i}\subset N_{2e}(C(\tau ))$, 
$\partial E_{i}\to \infty $, and
$|\langle E_{i},\omega \rangle| \neq 0.$ 
Here $\langle , \rangle$ denotes oriented intersection number. \\

\begin{pf}
By choosing $e > 2\delta$, we know that $N_{5e}(C(\tau)) \simeq \overset{\circ}{D} {}^{2} \times I$.\ 
Now, exhaust $\overset{\circ}{D} {}^{2} \times \{-1\}$ by concentric circles 
(say radius $r_i=1-\frac{1}{i+1}$, and call these curves $\tau_i$, and assume 
$w \cap { \overset{\circ}{D} {}^{2} \times \{-1\} }$ is in $\tau_1$ )
For any n, take very small neighborhood of $\tau_i$, $N(\tau_i)$, and change metric in 
very small neighborhood of $\partial U$, $N(\partial U)$, where $U=\tilde{M}-N(\tau_i)$, such that
$U$ has strictly convex boundary. Then by [MY], we have a least area disk in this metric, which is 
relatively least area disk in the original metric ,say $E_i$. Then by previous lemma, 
$E_i\subset N_e (\partial E_i)$, so $E_i\subset N_{7e}(C(\tau_i))$.
\end{pf}

\item

There exists a sequence $\{D_i\}$ of least area disks such that for all $n$, 
$D_i\subset N_{8e}(C(\tau))$, $\partial D_i\rightarrow \infty$ and 
$|\langle E_{i},\omega \rangle| \neq 0.$ \\

\begin{pf}
Since we did not change the metric outside very small neighborhood of $\tau_i$, we can cut down the 
size of $E_i$ such that $D_i\subset E_i$ are least area in $\tilde{M}$, and 
$\langle E_{i},\omega \rangle = \langle D_{i},\omega \rangle$.
\end{pf}

\item
After passing to a subsequence, $\{D_i\}$ converges to a lamination by least area planes which spans
$\tau$. \\

\begin{pf}

Let $\sigma $ be a $D^{2}$-limit lamination
obtained by applying Lemma 4.1 to $\{D_{i}\}$.   We still need to show that
each component of $\Si -\tau $ lies in a different complementary
region of $\sigma $.  If
$\omega _{1}\subset \BB ^{3}-\sigma $ is a properly embedded path 
connecting these two components, then since
$\omega _{1}\cap N_{2e}(C(\tau ))$  is compact and disjoint from
$\sigma $, it follows that for $i$ sufficiently large $D_{i}\cap \omega _{1}=\emptyset $.  
This contradicts the fact that for $i$
sufficiently large, $|\langle
\omega ,D_{i}\rangle| =|\langle \omega _{1},D_{i}\rangle|$.
\end{pf}

\item
if $\sigma$ spans $\tau\in \Si (\tilde M)$, then $\sigma\subset N_{9e} (C(\tau))$.\\

\begin{pf}

As we will prove in Lemma 5.2, for any i if 
$\partial D_i \subset K= N_\epsilon (\partial^- (N_7e(C(\tau))$
then $C(\partial D_i) \subset {N_{2\delta}(C(\tau) \cup A})$, where $A$ is union of 
geodesic segments from a to $a'=\pi(a)$, with $\pi:\partial D_i\rightarrow C(\tau)$ nearest 
point projection. But, since $\partial D_i\subset K$ then 
$A\subset N_{7e+\epsilon}(C(\tau))$. Then $C(\partial D_i)\subset N_{7e+\epsilon+2\delta}(C(\tau))$.
But $D_i\subset N_e(C(\partial D_i))$. So for any i,
$D_i\subset N_{9e}(C(\tau))$, assuming $e> \epsilon+2\delta$. 
as $\sigma=Lim\{D_i\}$, then $\sigma\subset N_{9e}(C(\tau))$.
\end{pf}

\end{enumerate}

\end{pf}

\section{ Genuine Laminations}

In second section, we get a $\pi_1$-invariant family of circles $\{C^+_x\}$
and $\{C^-_y\}$ in $S^2_{\infty }(\tilde{M})$. In third section, we spanned
these circles with lamination by least area planes in $\tilde M$. Now, we want to show
that these laminations indeed $\pi_1$-invariant, pairwise disjoint, and they
induce a pair of genuine laminations, $\Lambda^\pm$, on M.

\begin{thm}

There are laminations, $\hat{\sigma}^+ = \cup_{C^+_x}{\sigma}^+_x$ and 
$\hat{\sigma}^- = \cup_{C^-_y}{\sigma}^-_y$
in $\tilde{M}$ such that $\partial_\infty(\hat{\sigma}^+_x)=C^+_x$ and 
$\partial_\infty(\hat{\sigma}^-_y)=C^-_y$.
Moreover, these laminations are $\pi_1$-invariant, i.e. 
$\alpha(\hat{\sigma}^+_x) = \hat{\sigma}^+_{\alpha(x)}$, and 
$\hat{\sigma}^+_x \cap \hat{\sigma}^+_{x'}=\emptyset$.

\end{thm}

A short outline: First, we show that the lamination, for a fixed circle $\tau$ at
infinity, described in previous section does not intersect transversely with
the image of itself under a stabilizer of that circle. To show that, we use
the least area disks converging both laminations. There must be least area
disks in the sequences, which intersect transversely as the laminations. If
they intersect transversely, one of them must intersect the other one's
boundary. By fixing one of the discs, and choosing the other one very close
to the leaf of the lamination, we show that one of them cannot intersect the
other one's boundary.This is the first step. Then, we define the lamination
spanning the fixed circle $\tau$ as the union of the all the limiting laminations
of the sequences $\{\alpha_n\}$ in $\pi_1(M)$ such that $\alpha_n(\tau)\rightarrow \tau$.
By a similar method as above, we show that these images of the lamination
are pairwise disjoint. Then we can extend the lamination spanning a circle to whole family of circles 
by defining it the limit lamination for suitable sequence. Moreover, 
by construction they will be $\pi_1$-invariant.\\

\begin{pf}
Let $\tau \in \{C^+_x\}$, and
$G_\tau = Stab(\tau) = \{\alpha \in G_M \mid \alpha(\tau)= \tau \subset \Si \}$. 
We have a lamination by least area planes $\sigma_\tau$ by previous part, i.e. $\sigma_t$ is the limiting lamination
of sequence $\{P_i\}$, where $\partial P_i \subset \partial^- N_{5e}(C(\tau))$ and $\partial P_i \rightarrow \tau$ as
$i \rightarrow \infty$.

\begin{enumerate}

\item
$\sigma_\tau \cap \alpha(\sigma_\tau)=$ union of leaves of $\alpha(\sigma_\tau)$ and
$\sigma_\tau$, where $\alpha \in G_{\tau}$.\\

\begin{pf}
Assume in the contrary. Then there are leaves $L\in\sigma_\tau$ and $K \in \alpha(\sigma_\tau)$ such that
$L \cap K \neq \emptyset$ and the intersection is not the whole leave. So, it must be union 
of lines (maynot be disjoint), circles, and points. But, since $L$ and $K$ are least area planes then
the intersection cannot be a point, by maximum principle (Lemma 2.6 [HS]).
The intersection cannot be a circle, by exchange roundoff trick.

Now, we will prove it cannot be union of lines. By above discussion, 
we can find an intersection point 
$x$, where the intersection is transverse. 
By lemma 3.1., there are sequences of small disks 
$\{D_i\},\{E_i\}$ such that $D_i\in P_i$ and $E_i\in \alpha(P_i)=S_i$, 
$D_i\rightarrow D_x \subset L \cap B_\epsilon (x)$, $E_i\rightarrow E_x \subset K\cap B_\epsilon (x)$. 
Here, $\{P_i\}$ represents the least area disks defining $\sigma_\tau$.
Since $L$ and $K$ intersect transversely, for sufficiently 
large $i$ and $j$, $D_i$ and $E_j$ intersect transversely.

We claim that $\exists i_0, j_0$ such that $\forall i > i_0$, $D_i \cap E_{j_0} \neq \emptyset$.
Now, as $D_i\rightarrow D_x$, we can assume $d_H(D_i, D_x) \rightarrow 0$, where $d_H$ represents
Hausdorff distance. Since the intersection is transverse and $D_x$ and $E_x$ have bounded second 
fundamental form by [S], then $\exists {\epsilon'} << \epsilon$ such that
the distance between the sets $D_x - N_{\epsilon'} (D_x \cap E_x)$ and $E_x - N_{\epsilon'} (D_x \cap E_x)$
is greater than $\epsilon_1$, i.e. $E_x$ and $D_x$ does not get very close to each other away from the intersection.

Now, choose $i_0$ and $j_0$ such that $d_H(E_{j_0}, E_x)=\epsilon_2 << \epsilon_1$ and 
$d_H(D_i, D_x)< \epsilon_3 << (\epsilon_1- \epsilon_2)$. If $D_i$ does not intersect $E_{j_0}$ then 
$D_i$ belongs to a component of $B_\epsilon(x) - E_{j_0}$, but this contradicts to
$d_H(D_i, D_x)< \epsilon_3 << (\epsilon_1- \epsilon_2)$.

So,  we can assume that $\exists i_0, j_0$ such that 
$\forall i > i_0$, $P_i \cap S_{j_0} \neq \emptyset$. By the proof of the Lemma 3.3
$\partial D_i \subset N_\epsilon(\partial^-(N_7e(C(\tau))$ where $\partial^-$ represents the
lower part of the boundary. This $\epsilon$ comes from the process getting least area
disks from the relatively least area disks.

Now, choose sufficiently large $i > i_0$ such that $\partial P_i \cap S_{j_0}=\emptyset$ and 
$\partial P_i$ is very far from $\partial S_{j_0}$. If we show that 
$P_i \cap \partial S_{j_0} = \emptyset$, then this implies $P_i\cap S_{j_0}$ is not transverse,
as it is transverse one of them must intersect the other one's boundary. This will be a
contradiction and completes the proof of the claim.

\begin{lem}
There exist a uniform constant $C$ such that $P_i \cap T \subset N_C(\partial P_i)$ where
$T= N_\epsilon(\partial^-(N_{7e}(C(\tau))$.

\end{lem}

\begin{pf}
By lemma 3.2, we know that $P_i \subset N_e(C(\partial P_i)$. Now, consider $\partial P_i$, and
its nearest point projection to $C(\tau)$, say $\pi: \partial P_i \rightarrow C(\tau)$. 
Let $a'=\pi(a)\in C(\tau)$. Define $A= \bigcup_{a\in \partial P_i} \gamma_{aa'}$, where $\gamma_{aa'}$ 
represents the geodesic segment between $a$ and $a'$.

Now, we claim that $C(\partial P_i) \subset N_{2\delta}(C(\tau)\cup A)$. Let $x\in C(\partial P_i)$.
Then $\exists a,b \in \partial P_i$ such that $x\in \gamma_{ab}$. Now consider 
$a',b'\in C(\tau)$. Since $\tilde{M}$ is $\delta$-thin, 
$\gamma_{ab}\subset N_{\delta}(\gamma_{aa'}\cup \gamma_{a'b})$ and
$\gamma_{a'b}\subset N_{\delta}(\gamma_{a'b'}\cup \gamma_{bb'})$, so 
$\gamma_{ab}\subset N_{2\delta}(\gamma_{aa'}\cup \gamma_{a'b'}\cup \gamma_{bb'})$. 
But, $\gamma_{a'b'}\subset C(\tau)$ and $\gamma_{aa'}\cup \gamma_{bb'}\subset A$, this implies 
$C(\partial P_i) \subset N_{2\delta}(C(\tau)\cup A)$

Assuming $e> 2\delta$ , we can say that $N_e(C(\partial P_i)) \subset N_{2e}(C(\tau)\cup A)$.
Then $P_i \subset N_{2e}(C(\tau)\cup A)$. Consider $P_i\cap T$. Clearly, 
$T\cap N_{2e}(C(\tau)) = \emptyset$ as $7e-\epsilon >2e$. So if we prove
$N_{2e}(A)\cap T \subset N_C(\partial P_i)$, where C is independent of i, the claim follows.

Let $x\in N_{2e}(A)\cap T$. Then $\exists y \in \partial^-(N_{7e}(C(\tau)))$ such that $d(y,x)< \epsilon$
and $z\in \gamma_{aa'} \subset A$ with $d(y,z) < 2e+\epsilon$. Then $d(z,a) < 2e+2\epsilon$. Since
$d(y, C(\tau)) =7e$, $7e < d(y,a') \leq d(y,z)+d(z,a')$. Then $d(z,a') >5e-\epsilon$ and 
$d(a,a')\geq 7e-\epsilon$.

So, $d(a,x) < d(a,z) + d(z,y) +d(y,x) = 2e +2\epsilon + 2e + \epsilon + \epsilon = 4e + 4\epsilon =:C$
Then $P_i \cap T \subset N_C(\partial P_i)$. Lemma follows. 

\end{pf}

Now, we return to the proof of Step 1. Since $\alpha \in Stab(\tau)$ in $G_M$ acts as
isometry on $\tilde{M}$, $\alpha(\partial P_{j_0}) =\partial (\alpha(P_{j_0}))= \partial S_{j_0} \subset T$. 
Since $\partial P_i$ is very far away from $\partial S_{j_0}$
and $P_i\cap \partial S_{j_0} \subset P_i \cap T \subset N_C(\partial P_i)$, then
$P_i \cap \partial S_{j_0} = \emptyset$. Step 1 follows.

\end{pf}

Now, fix $\tau\in \{C^+_x\}$. 
Let $\sigma_0=\sigma$ as defined above. Define a set of sequences
$A:=\{ \{\alpha_n\} \subset \pi_1(M) | \alpha_n(\tau)\rightarrow \tau\}$.
Define $\sigma _{i+1}:=\bigcup_{\{\alpha_n\}\in A }lim \alpha_n(\sigma_i)$.
($lim \alpha_n(\sigma_i)$ is also lamination by least area planes, as we proved before.).
Then obviously,
$\sigma_0\subset\sigma_1\subset\sigma_2\subset .... \subset \sigma_n \subset ...$ with for any n 
$\partial_\infty\sigma_n=\tau$. Now, define $\hat\sigma_\tau=\sigma_\infty$ as described above.
Now, we will define the lamination for any circle $\tau'\in\{C^+\}$. By the construction of the 
lamination of $\lambda^\pm$ of $S^1$ [Ca], 
we know that the closure of the orbit of $\tau$ under the action of $\pi_1(M)$ on $\Si$ 
is the whole collection of circles $\{C^+\}$ ( Intuitively to get an idea what this means, consider a 
closed hyperbolic surface. Then take a nontrivial geodesic lamination on this surface. A dense leaf of this 
lamination lifts in universal cover ${\Bbb H}^2$ to an infinite geodesic. So the closure of the orbit of 
this leaf under $\pi_1(M)$ action will be the lift of whole lamination.)
So there exist a sequence 
$\{\alpha_n\}\subset\pi_1(M)$ such that
$\alpha_n(\tau)\rightarrow \tau'$. Then limit of the sequence $\alpha_n(\hat\sigma_\tau)$ will define 
another lamination $\hat\sigma_{\tau'}$ with
$\partial_\infty(\hat\sigma_{\tau'})=\tau'$.This is not very hard to see. Define a sequence of 
least area disks $\{S_n\}$ such that $S_n=\alpha_n(\hat\sigma_\tau)\cap N_c(C(\tau'))$.
Then these $S_n$'s will be sequences of least area disks whose boundaries are in $\partial^-(N_c(C(\tau')$.
Moreover, these sequence will converge to same lamination as the sequence $\alpha_n(\hat\sigma_\tau)$
since $\alpha_n(\tau)\rightarrow \tau'$.
On the other hand, this is independent of the choice of 
the sequence $\{\alpha_n\}$, by construction of $\sigma_\tau$.
So we define the family of laminations $\hat\sigma^+:=\{\hat\sigma_\tau | \tau\in\{C^+\}\}$.

In the following part,
we want to show that the union of the laminations $\{\hat{\sigma}^+_x\}$ constitutes a
lamination, $\hat{\sigma}^+$, in $\tilde{M}$\\

\item
Let $\mu, \omega \in \{ C^+ \} \subset \Si(\tilde{M})$.
Then $\hat{\sigma}_{\mu}\cap \hat{\sigma}_{\omega}= \emptyset$.\\

\begin{pf}
By Lemma 2.5, we know that that for any $C^+_x, C^+_{x'} \in \{C^+\} \subset \Si(\tilde{M})$,
the intersection is not transverse, i.e. $C^+_x\cap C^+_{x'}$ has only one component. 
if $\tau=\omega$ we are already done.
if not, the intersection is empty or at most one component. This means $C(\mu)$
and $C(\omega)$ cannot intersect transversely, one of them must lie one side of the other one.
Assume there are leaves of the $L\in \hat{\sigma}_{\mu}$ and $K \in \hat{\sigma}_{\omega}$,
intersecting transversely. We will adapt the proof of Claim 1.

First we modify the sequence of least area disks. As we defined above,
$\hat{\sigma}^+_\mu=\lim\alpha_n(\hat{\sigma}^+_\tau)$ and 
$\hat{\sigma}^+_\omega=\lim\beta_n(\hat{\sigma}^+_\tau)$  
where $\lim\alpha_n(\tau)=\mu$ and $\lim\beta_n(\tau)=\omega$. Consider the sequence 
$\{\alpha_n(\hat{\sigma}^+_\tau)\}$. Let $\{\overline{S_i}\}$ is the subsequence of $\alpha_n(\hat{\sigma}^+_\tau)$, where 
$\overline{S_i}$ is a least area plane in some $\alpha_n(\hat{\sigma}^+_\tau)$ and $\lim \overline{S_i}=\hat{\sigma}^+_\mu$. 
Now define a new sequence of disks, 
such that $S_i:=\overline{S_i}\cap N_{7e}(C(\mu))$. Since $\hat{\sigma}^+_\mu\subset N_e(C(\mu))$
$\lim \overline{S_i}=\lim S_i$. Similarly, if $\{\overline{P_i}\}$ is the subsequence of $\beta_n(\hat{\sigma}^+_\tau)$, where 
$\overline{P_i}$ is a least area plane in some $\beta_n(\hat{\sigma}^+_\tau)$ and $\lim \overline{P_i}=\hat{\sigma}^+_\omega$. Define
$P_i$ similarly.As $\hat{\sigma}_{\mu}$ and $\hat{\sigma}_{\omega}$ laminations by least
area planes, their intersection cannot be compact, i.e. they cannot intersect in a circle by exchange roundoff trick,
and they cannot intersect in a point by maximal principle for minimal surfaces. So the only possibility
the intersection must contain a line with endpoints $x,y\in I_{\mu\omega}$. Let's call this line
$l \subset K \cap L$ where K and L are least area planes in the laminations
$\hat{\sigma}_{\mu}$ and $\hat{\sigma}_{\omega}$ respectively.\\

\textbf{Case 1:} $\mu\cap\omega=\emptyset$.\

If $K \cap L \neq \emptyset$ then $K \cap L$ is a line, say $l$, by previous paragraph. But since $l=K \cap L$, then
$\partial_\infty(l) \subset \partial_\infty K \cap \partial_\infty L = \mu\cap\omega=\emptyset$, which is a contradiction.\\

%In this case, as we did before, we will concentrate on the converging disks around the transverse intersection of
%leaves $L\in \hat{\sigma}_{\mu}$ and $K \in \hat{\sigma}_{\omega})$. Now, choose sufficiently large N
%such that $S_N$ is sufficiently close to L and for all $n>N'$, $P_n$ and $S_N$ intersect transversely.
%We can find such N and N' as we did before. Now, since $\mu\cap\omega=\emptyset$, $C(\mu)\cap C(\omega)=\emptyset$
%as $\mu$ and $\omega$ are connected and live one side of each other (if there were intersecting geodesics, then WLOG
%$\mu$ must have 2 points in different components of complement of $\omega\Si$). Then by construction, as
%$\partial S_i \subset \partial^-(N_{7e}(C(\mu)$ and (WLOG assume negative side of $\partial(N_{7e}(C(\mu)$
%is the part far from $C(\omega)$, otherwise just choose the fixed disk
%from $\{P_i\}$ instead from $\{S_i\}$). By the proof of Lemma 4.2, we know that this converging disks
%$\{P_i\}$ are very close to the $C(\omega)$ ($P_i \subset N_{2e}(C(\omega)\cup A)$ where A is the annulus
%coming from union of geodesics between $\partial P_i$ and $C(\omega)$. So, if we choose n large enough then
%$P_n \cap \partial S_N$ must be empty. This means $S_N$ and $P_n$ do not intersect transversely. This is a contradiction.\\

\textbf{Case 2:} $\mu\cap\omega \neq \emptyset$.\

By Lemma 2.5, we know that if $\mu\cap\omega \neq \emptyset$, then the intersection has only one component,
say $\mu\cap\omega = I_{\mu\omega}$.
Now, we are at the only step which we use transverse orientability hypothesis. By transverse orientability,
the down sides  and up sides of the least area planes points the same sides as in Figure[4].

\begin{figure}
\mbox{\vbox{\epsfbox{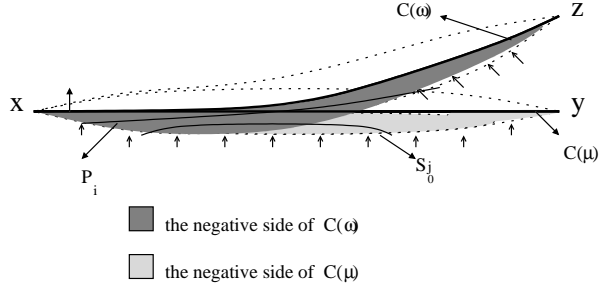}}}
\caption{\label{fig:figure2}
{2-dimensional picture of intersections of convexhulls of circles $\mu$ and $\omega$, which is represented in the figure
 by points $\{x,y\}$ and $\{x,z\}$, respecively. The line between x and y represents the convex hull of $\mu$, $C(\mu)$
 and the line between x and z represents the convex hull of $\omega$, $C(\omega)$. }}
\end{figure}

Now WLOG assume $\mu$ lies on the downside of $\omega$. 
Consider the sequences of least area disks converging to the 
transverse intersection, $P_i\rightarrow L$, and $S_j\rightarrow K$ as in the proof of Claim 1. Then 
again we can fix one disc, $S_{j_0}$ in one of the sequences 
and take another disc, $P_i$  intersecting the first one,very close to $L$ and the boundary of $P_i$
is very far from the $S_{j_0}$'s boundary. Remember by choice of the lamination,  
$\partial (S_j) \subset \partial^-(N_c(C(\mu)$ and $\partial(P_i) \subset \partial^-(N_c(C(\omega)$.
By Lemma 5.2. $P_i \cap T \subset N_c(\partial P_i)$. Then if we choose i sufficiently large $P_i$ cannot intersect 
$\partial S_{j_0} \subset T$. But this is a contradiction because if $P_i$ intersect $S_{j_0}$ transversely, $P_i$ 
must intersect $\partial S_{j_0}$
\end{pf}

\begin{figure}
\mbox{\vbox{\epsfbox{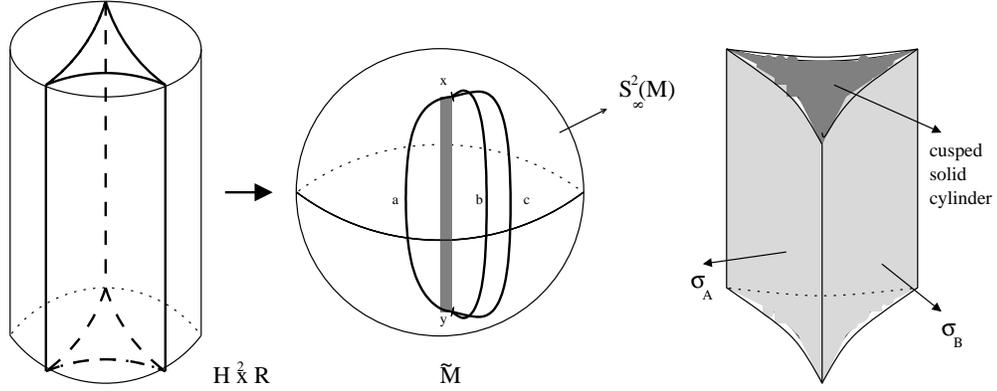}}}
\caption{\label{fig:figure3}
{$A,B,C \in \lambda^\pm$ 3 lines and they induce 3 circles in $\Si(\tilde{M})$, say $C_A,C_B,C_C$ where represent the circles through
points [x,a,y,b,x],[x,b,y,c,x],[x,c,y,a,x] respectively. When you span these circles at infinity with laminations $\sigma_A,\sigma_B,\sigma_C$
then there will be an infinite cusped solid cylinder, which is lift of cusped solid torus, between $\sigma_A,\sigma_B,\sigma_C$.}}
\end{figure}

\item
The lamination $\hat{\sigma}^+$ is $\pi_1$-invariant. i.e. for any $\alpha\in G_M$, 
$\alpha (\hat{\sigma}^+_\omega) = \hat{\sigma}^+_{\alpha(\omega)}$.\\

\begin{pf}
Let $\omega\in \{C^+_x\}$. Then by definition, $\hat{\sigma}^+_\omega=\lim\beta_n(\hat{\sigma}^+_\tau)$
and $\hat\sigma^+_{\alpha(\omega)}=\lim\gamma_n(\hat\sigma^+_\tau)$ where 
$\lim\beta_n(\tau)=\omega$ and $\lim\gamma_n(\tau)=\alpha(\omega)$. But, as we showed before, the definitions of 
$\hat{\sigma}^+_\omega$ and $\hat{\sigma}^+_{\alpha(\omega)}$ are independent of the choice of sequences, and clearly 
$\lim(\alpha(\beta_n))(\tau)=\alpha(\omega)$. This means 
$\hat\sigma^+_{\alpha(\omega)}=\lim(\alpha(\beta_n))(\hat\sigma^+_\tau)$, i.e.
$\alpha (\hat\sigma^+_\omega) = \hat{\sigma}^+_{\alpha(\omega)}$.
\end{pf}
\end{enumerate}
\end{pf}

So, by the $\pi_1$-invariance of the laminations, when we project down the
lamination via covering projection, we will get laminations $\Lambda^\pm$ in 
$M$. In other words, if $\pi:\tilde{M}\rightarrow M$ is covering projection,
then $\Lambda^\pm = \pi(\hat \sigma^\pm)$.

\begin{thm}
$\Lambda^\pm$ are a pair of transverse genuine laminations.

\end{thm}

\begin{pf}
First, we will prove $\Lambda^+$ is essential.
Each leaf $L_x$ of $\Lambda^+$ lifts 
to a 
surface 
$\tilde L_x$ in $\tilde{M}$ which is  a 
least  area {\it plane}, 
so $L_x$ is incompressible.  An end 
compression of $L_x$ would imply the existence  of a monogon in
$\tilde{M}$ connecting two very  close together subdisks of $\tilde L_x$ of 
very much larger area,
contradicting the fact that $\tilde L_x$ is least area as in  
Figure [3].  So, $\Lambda^+$ is essential.\\

Now, if we show that $\Lambda^+$ has gut regions, then we are done. If we look at
the lift of the lamination $\Lambda^+$, which is ${\hat\sigma}^+$, the lift of the complementary regions,
are the complementary regions of ${\hat\sigma}^+$. Consider that the family of circles $\{C^+_x\}$
are canonically coming from the lamination $\{\lambda^+_x\}$ in ${\hyp}^2$. By [Ca], 
there are some complementary regions which are ideal polygons in ${\hyp}^2$.
The image of the leaves in the boundary of this polygonal regions are going to be union of circles
such that one of them lies inside the other ones and each circle has at least 2 other circles with nontrivial
intersection.see figure [5].

Then the region between these circles will be asymptotic boundary of a complementary region.
Clearly, such a region cannot induce a interstitial bundle, so it must be gut region. 
So, $\Lambda^+$ is a genuine lamination.
 \end{pf}

\begin{rmk}
This additional hypothesis of transverse orientability is really necessary to work out this proof. 
It is because when you have 2 circles at infinity which intersects in an interval and their downsides and
upsides don't match up (i.e. the upside of one of them is the downside of the other one.), then the converging disks 
always intersects nontrivially no matter what happens, when there are least area planes in the laminations spanning
these 2 circles. So we cannot get a contradiction as above. See Figure[6].

\end{rmk}

\begin{figure}
\mbox{\vbox{\epsfbox{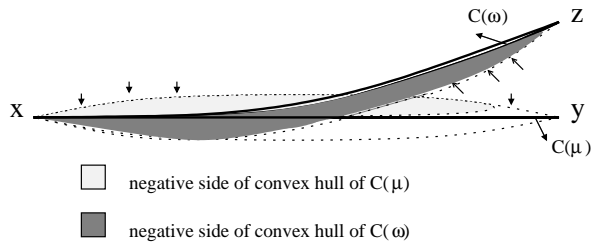}}}
\caption{\label{fig:figure4}
{ 2-dimensional picture of convex hulls of intersecting 2 circles at infinity whose negative sides don't match. }}
\end{figure}

\section{Topological Pseudo-Anosov Flows}

In this section we will show that by using the laminations defined in
previous section we could get a {\it Topological pseudo-Anosov flow} (TPAF)in the sense
of Mosher. \newline

In [Mo], Mosher defined TPAF and he proved
that if there is dynamic branced surface pairs in 3-manifold M , then we can
induce a TPAF. We will show the branched surfaces
carrying the laminations defined in previous section are actually a dynamic
pair,and by [Mo] we can induce a TPAF. The following definitions are from [Mo].

\begin{defn}

$\Phi$ is a TPAF if $\Phi$ has weak stable and unstable foliations,
singular along a collection of pseudohyperbolic orbits, and $\Phi$ has a Markov partition 
which is expansive in a certain sense (the latter condition is just relaxation of the expansive and 
contracting nature of smooth pseudo-Anosov flows.).

\end{defn}

This definition has two main purposes: First, it reflects many of the essential dynamic features
of a smooth pseudo-Anosov flow and and so many topological results about smooth 
pseudo-Anosov flows still hold. Second, it is much easier to verify in specific cases, like ours.

\begin{defn}

A {\it Dynamic Pair of Branched Surfaces} on a compact, closed 3-manifold M,
is a pair of branched surfaces $B^s, B^u\subset M$ in general position, disjoint from 
$\partial M$, together with a $C^0$ vector field V on M, so that the following conditions are satisfied.

\begin{enumerate}

\item
$(B^s,V)$ and $(B^u,V)$ are stable and unstable dynamic branched surfaces. (i.e. V is tangent to
$B^s$ and $B^u$ and along branch locus of $B^u$, $\Upsilon B^u$, points forward 
(from 2-sheeted side to 1-sheeted side and at crossing point 3-sheeted quadrant to 1-sheeted quadrant)
and along branch locus of $B^s$, $\Upsilon B^s$, points backward 
(from 1-sheeted side to 2-sheeted side and at crossing point 1-sheeted quadrant to 3-sheeted quadrant))

\item
V is smooth on M, except along $\Upsilon B^s$ where backward trajectories locally unique, and 
along $\Upsilon B^u$ where forward trajectories locally unique.

\item
Each component of $M - (B^s \cup B^u)$ is either a pinched tetrahedron or a solid torus. In solid
torus piece, V is circular. See Fig[??]

\item
Each component of $B^u - B^s$ and $B^s - B^u$ is  an annulus with cusped tongues, see figure[??].
On components of $B^u - B^s$,the annulus is a sink for V (all forward trajectories of V after a time is 
in the annulus.), and similarly on components of $B^s - B^u$,the annulus is a source for V (all backward trajectories of V after a time is 
in the annulus.).

\item
No two solid torus components of $M - (B^s \cup B^u)$ are glued to each other, i.e. 
the closures of solid torus components are disjoint.

\end{enumerate}

\end{defn}

Now, let $\Lambda ^{\pm }$ be the genuine laminations defined in previous
section. Let $B^{\pm }$ be the branched surfaces carrying $\Lambda ^{\pm }$.
We want to show that $B^{\pm }$ are dynamic pair. Here, $B^+$ and $B^-$ correspond to $B^s$
and $B^u$, respectively.

\begin{lem}

$\Lambda^\pm$ are very full laminations in M, i.e. gut regions are solid tori.

\end{lem}

\begin{pf}
This is true as the gut regions are coming from the ideal polygons of the lamination
$\lambda^\pm \subset S^1_infty({\Bbb H}^2$. These ideal polygons induces circles at infinity as in the Figure[5].
So the gut regions are the region between the lamination spanning this circles. On the other hand for each ideal
polygon we have an element $\alpha$ in $\pi_1(M)$ fixes this ideal polygon (the topological pseudo-Anosov elements in [Ca]).
Then $\alpha$ fixes the two common points of all the circles coming from the each side of ideal polygon. So, the gut region must be 
a solid tori whose core is homotopic to the element $\alpha$. So, the gut regions are solid tori.
\end{pf}

Now, recall that the lamination $\Lambda^\pm$ is coming from universal cover
and the lifting laminations $\tilde{\Lambda^\pm}$ are laminations by least
area planes. Let P be a least area plane in the lamination and let
$\partial_\infty(P)= \tau \in \{C^+_x\}$. Then we have special point $a\in
\tau\subset \Si(\tilde{M})$. By Lemma 2.4, $\tau=
q \circ p(\partial_\infty(l^+_x \times I)))$ and by proof we know that
$q \circ p(\partial(l^+_x)\times I \cup l^+_x \times \{\infty\})$is a point and we
define this point as special point in $\tau$.

%Now, consider $P\subset\tilde{M}$ with the induced metric. Define a geodesic foliation of P by geodesic 
%whose one endpoint is the special point $a\in\tau=\partial_infty(P)$ and the other endpoint is in $\tau-\{a\}$.
%Induce a vector field from this foliationsuch that each vector points towards the special point on line. 
%Then do same thing for any least area plane in $\Lambda^+$. Now we have a vector field on $\Lambda^+$.
%Now, define a similar vector field on $\Lambda^-$ but this time the vectors points out from the special point on line. 

Let $B^\pm$ branched surfaces carrying the genuine laminations $\Lambda^\pm$ such that branch locus of 
$B^\pm$ is transverse to the $B^+ \cap B^-$.

\begin{thm}

If $B^\pm$ branched surfaces carrying the genuine laminations $\Lambda^\pm$ then $B^\pm$ are a 
dynamic pair of branched surfaces. So, there is a topological pseudo-Anosov flow on M by [Mo]. 

\end{thm}

\begin{pf} There are 5 steps.\\

\begin{enumerate}

\item
(Structure of $B^\pm$ and M)\\
$M-B^+$ is union of cusped tori, $Q^+_i$. Similarly, $M-B^-$ is also union of cusped tori, $Q^-_i$. Moreover
$Q^+_i \cap Q^-_i = T_i$ where $T_i$ represents solid torus gut region for $\Lambda^\pm$. $M-(B^+\cup B^-)$
consists of solid tori ( not cusped) and pinched tetrahedron Figure[7].On the other hand,
components of $B^+ - B^-$ and $B^- - B^+$ are annuli with "cusped" tongues as in the figure [6].\\

\begin{figure}
\mbox{\vbox{\epsfbox{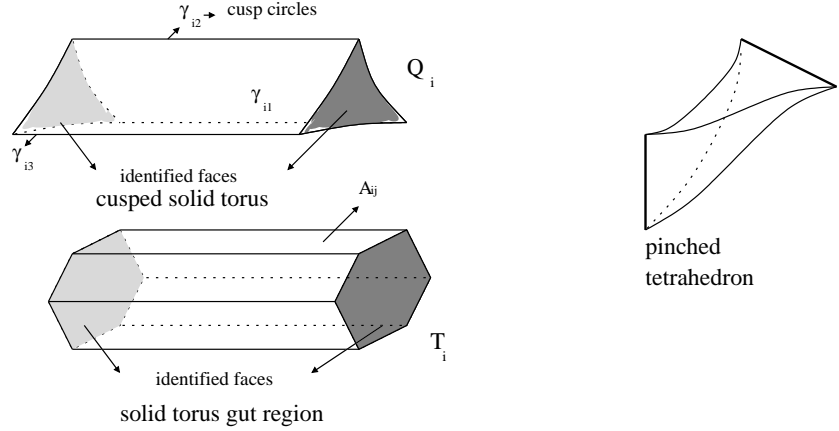}}}
\caption{\label{fig:figure5}
{Shapes of cusped solid torus and pinched tetrahedron pieces in $M-(B^+ \cup B^-)$. Solid torus gut region is
intersection of 2 transverse cusped solid torus pieces.}}
\end{figure}

Since $\Lambda^\pm$ are very full laminations, then $M-B^+= \bigcup_{i=1}^n Q_i^+$ where $Q_i^+$ 
represents cusped solid torus piece, see Figure [7].
similarly $M-B^- = \bigcup_{i=1}^n Q_i^-$. Moreover for any i, $Q_i^+ \cap Q_i^- = T_i$ where $T_i$ is the
(noncusped) solid torus gut piece of the lamination. As we have seen above, these gut regions, $T_i$,
comes from $r_i$ sided ideal polygons in $\lambda^\pm$, as we call them $r_i-prong$. Then these cusped torus pieces,
$Q_i^\pm$ have $r_i$ cusp circles, say $\gamma^\pm_{ij}$ $1 \leq j \leq r_i$. In the boundary of corresponding gut
region $T_i$, there are $2r_i$ parallel circles, coming from the intersection $Q_i^+\cap Q_i^- =T_i$. These
$2r_i$ circles in the boundary of solid tori $T_i$, bounds $2r_i$ annuli in $\partial(T_i)$ and these annuli
alternatingly in $B^+$ and $B^-$. if it is in $B^+$, we will call them $+annulus$ and if it is in $B^-$
then we will call them $-annulus$. 

\begin{figure}
\mbox{\vbox{\epsfbox{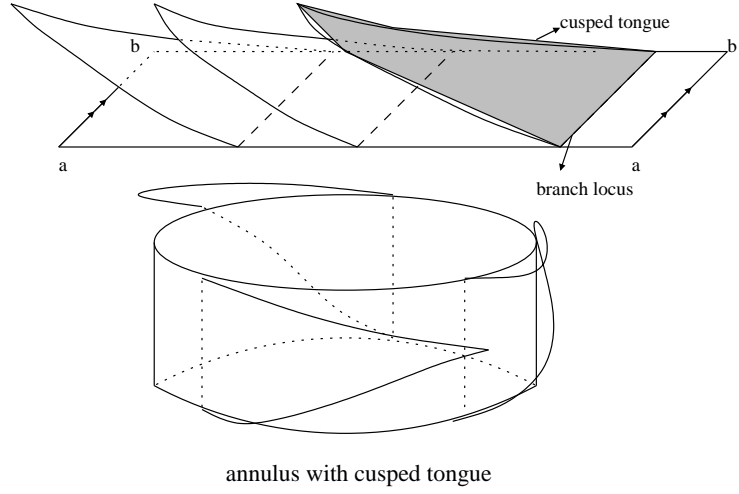}}}
\caption{\label{fig:figure6}
{Shape of annulus with 3 cusped tongues.}}
\end{figure}

Take a $+annulus$ in $\partial(T_i)$. This annulus comes from the intersection 
of a cusp in $Q^-_i$ and $B^+$. So we can index these annuli, by just considering the indexing of cusps coming
from $\gamma^\pm_{ij}$. So for each $+annulus$ there is a $\gamma^-_{ij}$ and for each $-annulus$, there is a 
$\gamma^+_{ij}$. Then call the $+annuli$ corresponding to $\gamma^-_{ij}$ as $A^+_{ij}$ and similarly define 
$A^-_{ij}$. Now, we have $\partial(T_i)= \cup_{j=1}^{r_i}A^+_{ij} \bigcup \cup_{j=1}^{r_i}A^-_{ij}$.

Each cusp circle $\gamma^-_{ij}$ and $A^+_{ij}$ defines a cusp, say $C^-_{ij}$, in $Q^-_{ij}$ and 
similarly $C^+_{ij}$, in $Q^+_{ij}$. Then the cusped solid torus
$Q^-_i = T_i \bigcup \cup_{j=1}^{r_i}C^-_{ij}$ and similarly $Q^+_i = T_i \bigcup \cup_{j=1}^{r_i}C^+_{ij}$
See Figure[9].

\begin{figure}
\mbox{\vbox{\epsfbox{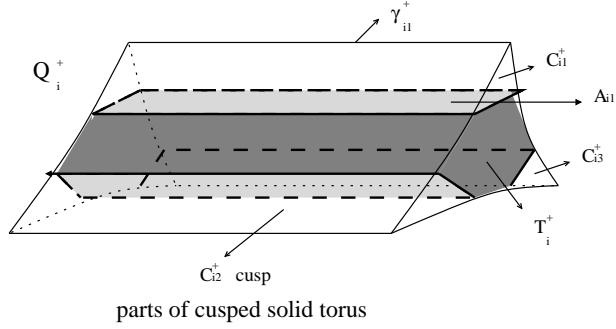}}}
\caption{\label{fig:figure7}
{$Q_i^+=T_i \cup (\bigcup_{j=1}^{3} C_{ij}^+)$ }}
\end{figure}

Now, let's describe the pieces of $B^+-B^-$. we claim that these pieces are annuli with "cusped" tongues as in the 
figure [8]. Consider $M-B^+= \bigcup_{i=1}^n Q_i^+$. then  $\bigcup_{i=1}^n \partial(Q_i^+) \supseteq B^+$. So if we understand,
how +cusped tori and -cusped tori intersect, then we can easily decribe the components of $B^+-B^-$. But as we
mentioned above, these intersections produce solid tori gut regions and cusps. This means that
components of $B^+-B^-$ will have one of annulus $A^+_{ij}$ and the remaining part of the component 
will be in the cusp $C^+_{ij}$. It is easy to see that these parts in the cusp will be the cusped 
tongues coming from the other sections of the branched surface $B^+$ as in the
Figure [8] (section of a branched surface is the components of branched surface - branch loci,
$B^+ - \Upsilon B^+ $).

The other claim is that the components of $M- (B^+ \cup B^-)$ are solid tori and pinched tetrahedra. 
Consider the following trivial set theoretic equivalences.
$M-(B^+ \cup B^-)= (M-B^+)\bigcap(M-B^-)= (\bigcup_{i=1}^n Q_i^+)\bigcap (\bigcup_{i=1}^n Q_i^-) =
(\bigcup_{i=1}^n (Q_i^+ \cap Q_i^-)) \bigcup (\bigcup_{i \neq k} (Q_i^+ \cap Q_k^-))
= (\bigcup_{i=1}^n T_i) \bigcup (\bigcup_{i \neq  k} (Q_i^+ - T_i) \cap (Q_k^- - T_k))
=(\bigcup_{i=1}^n T_i) \bigcup (\bigcup_{ij \neq  kl} C_{ij}^+ \cap C_{kl}^-)$

Now, the first part of the union comes from the equality $Q_i^+\cap Q_i^- =T_i$, intersection of cusped solid tori with same indices
is the corresponding solid torus gut region. In the latter part of the union we just used the definitions in the first paragraph:
$Q_i^+ - T_i = \bigcup_{j=1}^{r_i} C_{ij}^+$, the cusped solid tori are the union of solid tori gut regions and the cusps.

So, if we can understand $C_{ij}^+ \cap C_{kl}^-$ for i and k different, then we will finish this step.
we claim that this intersections give us the pinched tetrahedra components. Consider the Figure [10].
As it can be seen there the intersection of the cusps of different cusped solid tori is in general position (by assumption,
$\tau=B^+ \cap B^-$ is transverse to the branch loci of the branched surfaces, $\Upsilon B^\pm$.). Fix a cusp
$C_{ij}^+$ in $Q_i^+$. Now, consider the intersection of $C_{ij}^+$ with the other regions. Obviously, since
this region lives already in the complement of $B^+$, $\bigcup_{i=1}^n Q_i^+$ , no region in the complement of $B^+$ intersect $C_{ij}^+$.
Now, consider the intersection with $\bigcup_{i=1}^n Q_i^-$. Since solid tori gut regions are disjoint from cusps
then only cusps of the negative cusped solid tori will intersect our region $C_{ij}^+$. 

\begin{figure}
\mbox{\vbox{\epsfbox{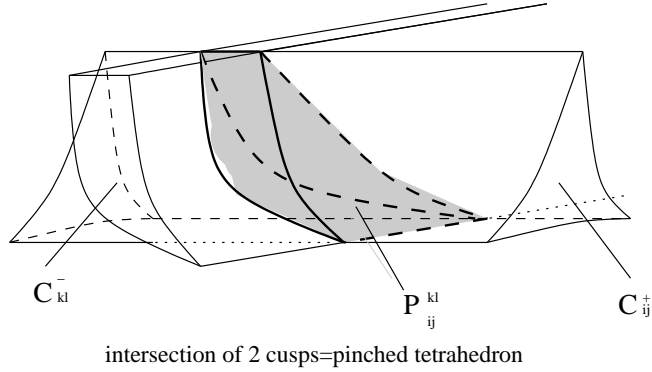}}}
\caption{\label{fig:figure8}
{Intersection of 2 cusps, $C_{ij}^+ \cap C_{kl}^-$, is a pinched tetrahedron, $P_{ij}^{kl}$. }}
\end{figure}

Recall that the cusps are topologically just a cusped (in one vertex) triangle $ \times S^1$. the cusp vertex $ \times S^1$
corresponds cusp circle which is in branch locus of $B^+$, $\Upsilon B^+$,
and the opposite side of triangle $ \times S^1$ corresponds the annulus in $B^-$. Now the negative cusps intersect our cusp circle in intervals
and the annulus have some interval parts of branch locus of $B^-$. These intervals will constitute the cusped sides of a tetrahedra intersections,
and the intersections of positive and negative cusps will be pinched tetrahedra. So, the components of $M- (B^+ \cup B^-)$ are solid tori 
and pinched tetrahedra as claimed.\\

\item
We can define vector field X on M which is tangent to $\tau=B^+ \cap B^-$ and $B^+$ and $B^-$.\\

First, we will define the vector field on train track $\tau=B^+ \cap B^-$ and then we
will extend first to $B^+ - B^-$ and $B^- - B^+$ naturally.\\

\begin{itemize}

\item
X on $\tau$:\\

It is not obvious that we can define a vector field on a train track, see Figure[11].

\begin{figure}
\mbox{\vbox{\epsfbox{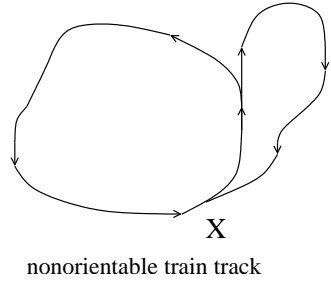}}}
\caption{\label{fig:figure9}
{we cannot define a vector field on this train track.}}
\end{figure}

This is indeed same thing with orienting each segment in train track consistently.
First we will show that we can define canonically a vector field on $\tau$ by using the circles at infinity in 
universal cover. If we consider the lift of branched surfaces in universal cover $\tilde {B^\pm}$, we can see the the
intersection train track lifts to infinite lines asymptotic to the end of lifts of solid tori, which are the special points (defined above) of
corresponding circles at infinity, i.e. each infinite line limits to one positive special point (special point in a positive circle 
$C^+_x$ at infinity)
and to one negative special point ( to see intuitively consider the quasi-isometric picture of $\tilde{M}$ as 
$\hyp^2 \times \BR$, and the infinite lines starts from bottom disk and ends in top disk)
So clearly we can orient each infinite line from a negative special point to positive special point. 
Now, we will induce consistent orientation of each segment of $\tau$ using these orientation of lines in
$\tilde{\tau}$. Take a line segment $I \subset \tau$ and consider a lift of this line segment 
$\tilde{I} \subset \tilde{\tau}$ in universal cover. Clearly, we can orient the circles in $\tau$ 
which are in boundary of solid torus gut regions (for each $T_i$, there are $2r_i$ circles in 
$\partial(T_i)$ which are also in $\tau$) parallel to the core of the gut region. 

Now the only remaining 
part of $\tau$ to orient is the line segments connecting these circles. Consider the the quasi-isometric 
picture of $\tilde{M}$ as $\hyp^2 \times \BR$. In this picture as we have seen in Section 2, 
the family of circles at infinity $\{C^\pm\}$, comes from $\partial_\infty(\lambda^\pm \times \BR)$
by collapsing $\lambda^+$ in $\hyp^2 \times \{+\infty\}$ and by collapsing $\lambda^-$ in 
$\hyp^2 \times \{-\infty\}$.Since $B^\pm$ carries the laminations $\Lambda^\pm$ 
(i.e. $\Lambda^\pm \subset N_\epsilon(B^\pm)$),$\partial_\infty(\Lambda^\pm)=\partial_\infty(B^\pm))$.
So, if you take two "close" leaves of lifts of $\tilde{B}^+$ they will intersect in an interval not 
containing their special point of both circles and they will start to differ from their special point
(Recall that every circle at infinity, $\{C^\pm\}$, has a special point which is the image of the 
endpoint of corresponding leave of $\lambda^\pm$) See Figure[12]. This is true for $\tilde{B}^-$ as well. So,
for the circles corresponding to the sides 
of ideal polygons in $\lambda^\pm$ and corresponding circle at infinity of the leaves in $\tilde{B}^\pm$
containing boundaries of solid tori gut regions, they have both negative and positive special points, and
as in previous paragraph we oriented the core of solid torus as from negative special point to 
positive special point. See Figure[13]

\begin{figure}
\mbox{\vbox{\epsfbox{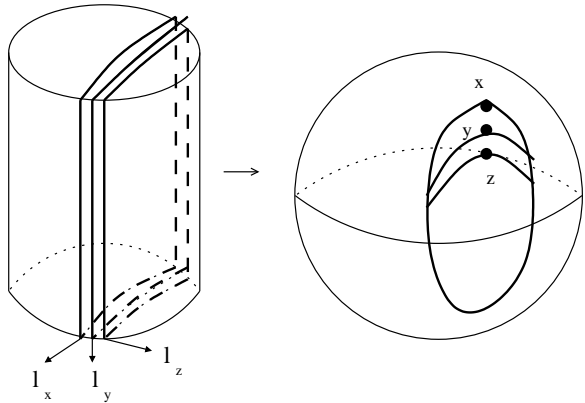}}}
\caption{\label{fig:figure10}
{Induced circles from 3 generic leaves,  $l_x,l_y,l_z \in \lambda^+$ (not a boundary of ideal polygon in the complement of $\lambda^+$)
in $\Si(\tilde{M})$
  }}
\end{figure}

\begin{figure}
\mbox{\vbox{\epsfbox{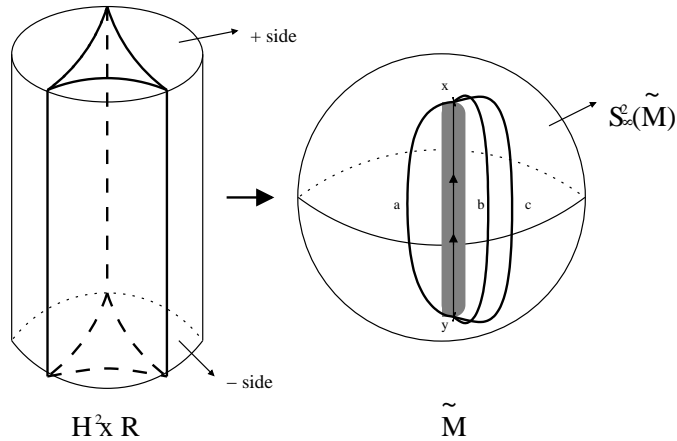}}}
\caption{\label{fig:figure11}
{Induced circles from 3 nongeneric leaves,  $l_a,l_b,l_c \in \lambda^+$ (sides of ideal polygon in the complement of $\lambda^+$)
in $\Si(\tilde{M})$
 }}
\end{figure}

Now, observe that in $\hyp^2 \times \BR$ picture, the lift of branch locus, $\tilde{\Upsilon}^+$ (which are lines
as loops in branch locus are essential), in $\tilde{B}^+$ branches towards positive side of 
$\hyp^2 \times \BR$, and similarly, $\tilde{\Upsilon}^-$ in $\tilde{B}^-$ branches towards negative side of 
$\hyp^2 \times \BR$, see figure [14]. This is very easy to see if the laminations are geodesic planes in $\hyp^2 \times \BR$, because of the tightness.
But in our situtation the tightness comes from being least area planes, which works in our situation as well.
In other words, we know that the close circles at infinity, say $C^+_1,C^+_2$ starts to diverge from each other from their 
special points and this will cause inside $\tilde{M}$ the leaves $L_1, L_2$ of lamination $\Lambda^+$ 
will be close to each other for some time but they will start to diverge from each other after a lift of
intersititial annulus. See figure [15]. On the other hand this intersititial annulus corresponds in branched surface
literature a branch locus. Now, we want to say that this branchings towards upside for $B^+$ and towards downside for
for $B^-$. This is true as at infinity diverging starts at positive side and inside we have tightness coming from
the lamination being by least area planes. 

Now, let's come back to $\tau$. For a line segment in $\tau$ starts from $\Upsilon B^+$ and ends in
$\Upsilon B^-$ will be as in Figure[16] . So we will orient this line segment from $\Upsilon B^+$ to
$\Upsilon B^-$. Then our quasi-isometric picture of $\tilde{M}$ as $\hyp^2 \times \BR$ shows that 
the orientation on each line of $\tilde{\tau}$ is coherent, and when we project 
it to the original manifold, we can easily get a vector field on our train track $\tau$.\\

\begin{figure}
\mbox{\vbox{\epsfbox{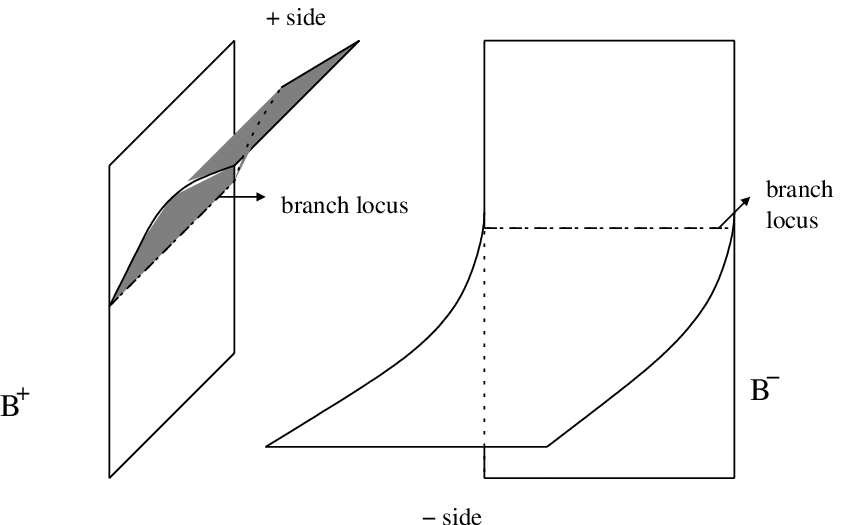}}}
\caption{\label{fig:figure12}
{Shape of neighborhood of $\tilde{\Upsilon B^\pm}$ in $\tilde{B^\pm}$ in $\hyp^2\times\BR$ picture of $\tilde{M}$. }}
\end{figure}

\begin{figure}
\mbox{\vbox{\epsfbox{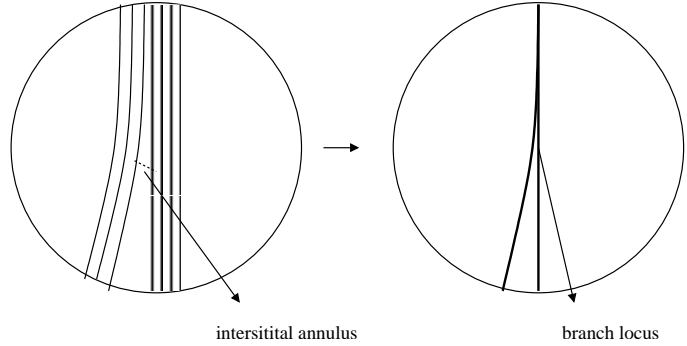}}}
\caption{\label{fig:figure13}
{2 dimensional picture of laminations and branched surfaces carrying them. Intersititial annulus becomes branch locus.}}
\end{figure}

\item
Extending X to the components of $B^+ - B^-$ and $B^- - B^+$:\\

By the first step we know that the components are annuli with cusped tongues. Now fix a component. 
Then its boundary will be in $\tau$, and we already defined X on $\tau$. Now, as we pointed before, 
since we induced X on $\tau$ from universal cover's boundary at infinity, there is no consistency problem.
i.e. since X is well-defined on $\tau$, on the boundary of annulus of component, they must be parallel, and 
on boundary of cusped tongues they are consistent. So we can easily extend first on annulus such that
each integral integral curve of X on annulus is closed as in boundary
(as X is parallel on two circles of the boundary),
and then on cusped tongues. If we have a +annulus with cusped tongue then X on $\tau$ points away from the ideal vertex towards
the annulus, and we can extend X to the cusped tongue with integral curves starting at ideal vertex, tangent to the sides contatining 
ideal vertex, and ending in the opposite side of ideal vertex, which is a segment of $\Upsilon B^-$. Similarly,
we can extend X to -annulus with cusped tongue.\\

\item
Extending X to whole manifold by defining on the solid torus and pinched tetrahedron pieces.\\

We have defined X on whole $B^\pm$. As we proved before components of $M- (B^+\cup B^-)$ are solid tori and 
pinched tetrahedra. First, let's extend X to pinched tetrahedron pieces.
Fix a pinched tetrahedron P. $\partial P$ consists of 4 cusped tongues, one couple comes from a 
positive annulus with cusped tongues (the component is in $B^+-B^-$)  and the other couple comes from 
negative annulus with cusped tongues(the component is in $B^- -B^+$). Now, there are 2 cusped segments in P,
one is an interval $I^+$ in $\Upsilon B^+$, and the other is an interval $I^-$ in $\Upsilon B^-$. Now, by our definition of X on 
$\tau$, and it's canonical extension to the components of $B^+ - B^-$ and $B^- - B^+$, X points inside to P
on $I^+$ and points outside from P on $I^-$. Then, it is clear that we can extend X to whole P such that,
X will be tangent to $\partial P$ and any integral curve of X in P starts from $I^+$ and ends in $I^-$.

\begin{figure}
\mbox{\vbox{\epsfbox{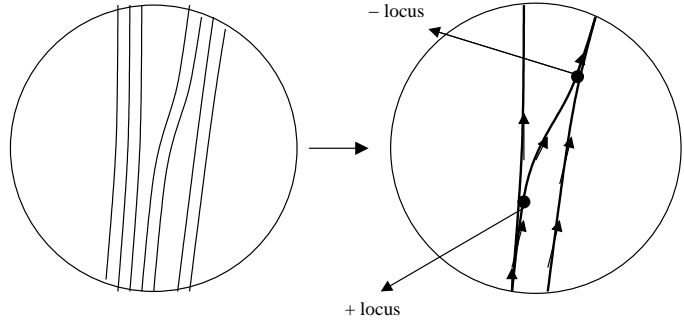}}}
\caption{\label{fig:figure14}
{Orienting the train track $\tilde{\tau}$ }}
\end{figure}

Now, fix a solid torus $T_i$ in $M - (B^+\cup B^-)$. As above, $\partial T_i$ consists of $2r_i$ annuli
from $B^\pm$. Boundaries of these annuli are $2r_i$ closed curves in $\tau$, and the definition of X on 
these annuli canonically comes from the definition of X on these circles. But, we defined X on $\tau$
by using the lift of $T_i$ to universal cover, and on each of these closed curves on $\partial T_i$ X is parallel
to the orientation of the core curve of $T_i$. So on each annuli the integral curves of X are closed and have 
same orientation with the core curve of $T_i$. It is obvious that we can simply extend X to $T_i$ such that 
each integral curve is closed and oriented parallel to core curve (i.e. the integral curves on solid
torus $T_i$ will be the trivial one dimensional foliation.).

Now, we have to check that X is continuous on M, i.e. there is no consistency problem with 
the definition of X on different components. Since there cannot be any problem inside the pinched 
tetrahedron and solid torus pieces, we should check only the boundaries of these pieces which are 
$B^+ \cup B^-$. But already we have induced X from the boundaries of the pieces, X is also continuous 
on the boundaries, i.e. $B^+ \cup B^-$. So, X is a $C^0$ vector field on M and it is tangent to 
$B^+\cup B^-$, such that X points inside to $B^-$ on $\Upsilon B^-$ and points outside from $B^+$
on $\Upsilon B^+$.\\

\end{itemize}

\item
There is no face gluings between solid torus gut regions, $T_i$, i.e. torus pieces of $M-(B^+\cup B^-)$ are separated.\\

Assume there is a face gluing between two solid torus components, say $T_i, T_j$. This means
there is a common annulus piece in $\partial(T_i) \cap \partial(T_j)$. When we look at the lifts of
$T_i$ and $T_j$ to the universal cover, we see that there is only one plane component of the lift 
of $B^+$ or $B^-$ separating these two lifts $\tilde{T_i}$ and $\tilde{T_j}$. On the other hand, that means
the boundary at infinity of this plane is isolated in both sides. This is not hard see, as these solid tori
components comes from ideal polygons in the lamination of circle $\lambda^\pm$. See figure[17].

\begin{figure}
\mbox{\vbox{\epsfbox{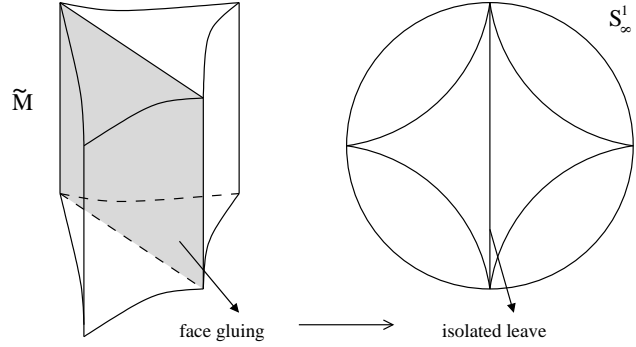}}}
\caption{\label{fig:figure15}
{Face gluings implies isolated leaves. Left ideal triangle of $\lambda^+ \subset S^1_\infty$ induce one cusped solid torus, and right
ideal triangle induce the other cusped solid torus.}}
\end{figure}

But this is contradiction since isolated circle at the boundary at infinity means isolated leaf of
the lamination $\lambda^\pm$ and we already know by [Ca] that $\lambda^\pm$ has no isolated leaves.\\

\item
$B^\pm$ are dynamic pair of branched surfaces.\\

The steps 1, 2, 3 proves the first 5 conditions of dynamic pair of branched surfaces and the step 4 shows
the last condition of dynamic pair of branched surfaces. So, $B^\pm$ are dynamic pair of branched surfaces.

\end{enumerate}

\end{pf}

This means if M is an atoroidal 3-manifold admitting uniform 1-cochain, then there is a 
TPAF on M induced by the uniform 1-cochain. If we consider uniform 1-cochains as generalization of sliterings
this is a generalization of a theorem of Thurston [Th]: if an atoroidal 3-manifold M slithers around circle 
then there is a pseudo-Anosov flow on M, transverse to the uniform foliationinduced by slithering.
In our setup, the uniform foliation corresponds the coarse foliation of $\tilde M$ induced by uniform 1-cochain.

\section{Concluding Remarks}

The transverse orientability condition on uniform 1-cochain is a little bit strong and disturbing.
To get rid of this hypothesis, one can try different approaches. One of them could be the below conjecture.\\

\textbf{Conjecture:}
Let M be Gromov hyperbolic 3-manifold, and $\alpha$ and $\beta$ are two simple closed curves
in $\Si(\tilde{M})$. If the least area planes K, and L spanning $\alpha$ and $\beta$, respectively, intersect
transversely in a line l which limits $\{x,y\} \subset \Si(\tilde{M})$, then the circles $\alpha$ and $\beta$
intersect transversely at $\{x,y\}$.\\

This might seem a very optimistic conjecture because in one less dimension this is not true, as geodesics may
intersect and stay in bounded Hausdorff distance in Gromov hyperbolic manifolds. But,
2-dimensionality of the objects might be very crucial and essential here. 
If this conjecture was true, the above theorem would follow easily as the planes in laminations would automatically be
pairwise disjoint. Moreover, this conjecture would make this technique so powerful that
to get an essential lamination in Gromov hyperbolic manifolds would be equivalent to get a $\pi_1(M)$ invariant
family of circles at infinity.

On the other hand, the minimal surface techniques and results in this paper are indeed original in the sense that
it starts with an algebraic condition on fundamental group $\pi_1(M)$, like admitting a function to $\BR$, uniform 1-cochain,
and ends up with two real topological object in the manifold M, like genuine laminations and topological pseudo-Anosov flow.
Of course, most of the work has been done by Calegari in his beautiful paper [Ca].

In last five years, we have seen three breakthrough results of nonexistence of some promising structures in 3-manifolds. Roberts, Shareshian, and Stein
proved that there are hyperbolic manifolds without taut foliations, [RSS]. By that time, it was believed that taut foliations are very abundant in 3-manifolds,
it might even be enough for weak hyperbolization. By [RSS], we saw that this is not true. The next promising structure for weak hyperbolization was essential
laminations. Calegari and Dunfield showed that tight essential laminations in atororidal manifolds induce circle action of the fundamental group and 
the fundamental group of the Weeks manifold does not act on circle. So this is the first example of hyperbolic manifolds without tight essential
laminations. Finally, Fenley showed that there are hyperbolic manifolds without any essential laminations, [Fe]. Taut foliations and essential laminations
were expected to provide a positive answer for weak hyperbolization before these results. 

Similarly, after Thurston's paper on slitherings, [Th], then their generalization as uniform 1-cochains by Calegari, and
abundance of bounded 1-cochains by geometric group theory, uniform 1-cochains might also be considered
as a promising tool for weak hyperbolization. The above paper of Calegari and Dunfield also show that there are hyperbolic manifolds
without uniform 1-cochains. Since uniform 1-cochains on atoroidal manifolds induce faithful circle action of fundamental group by [Ca],
they showed that the fundamental group of the Weeks manifold does not act on circle, so Weeks manifold cannot admit uniform 1-cochain.

When we started this problem, [CD] and [Fe] were not published yet, and we believed that by proving these results, we can contribute to Thurston's 
and Calegari's promising program for weak hyperbolization. After [CD] and [Fe], one can look at our results as another way of proving
nonexistence of uniform 1-cochains in some hyperbolic manifolds, up to transverse orientability condition. 
This is because by our work transversely orientable uniform 1-cochains induce genuine
laminations and by [Fe], there are hyperbolic manifolds without genuine laminations.

\end{document}